\input ./style/arxiv-general.cfg
\documentclass[MSNbibl,number,citesort,seceqn,dvips]{arxbj}
\makeatletter
   \@ifpackageloaded{graphicx}{}{\usepackage{graphicx}}
\makeatother

\volume{22}
\issue{4}
\pubyear{2016}
\firstpage{2548}
\lastpage{2578}
\doi{10.3150/15-BEJ737}
\docsubty{FLA}

\makeatletter
\def\vfrac#1#2{(#1)/#2}
\def\sfrac#1#2{#1/#2}
\newcommand{\xrightarrow}[1]{\stackrel{\mathrm{#1}}{\to}}

\newcommand{\rrvert}{\vert}

\newcommand{\rrVert}{\Vert}
\newcommand{\llvert}{\vert}
\newcommand{\llVert}{\Vert}

\newtheorem{thmo}{Theorem}[section]
\newtheorem{lem}[thmo]{Lemma}
\newtheorem{prop}[thmo]{Proposition}
\newproclaim{defn}[thmo]{Definition}
\newremark{example}[thmo]{Example}
\newremark{notn}[thmo]{Notation}
\makeatother

\begin{document}
\begin{frontmatter}

\title{An $\alpha$-stable limit theorem under sublinear expectation}
\runtitle{An $\alpha$-stable limit theorem under sublinear expectation}

\begin{aug}
\author[A]{\inits{E.}\fnms{Erhan}~\snm{Bayraktar}\corref{}\thanksref{e1}\ead[label=e1,mark]{erhan@umich.edu}}
\and
\author[A]{\inits{A.}\fnms{Alexander}~\snm{Munk}\thanksref{e2}\ead[label=e2,mark]{amunk@umich.edu}}
\address[A]{Department of Mathematics, University of Michigan, Ann Arbor,
MI 48109, USA.\\ \printead{e1,e2}}
\end{aug}

%
\received{\smonth{1} \syear{2015}}
%
\revised{\smonth{5} \syear{2015}}

%
\begin{abstract}
For $\alpha\in ( 1 , 2  )$, we present a generalized
central limit theorem for $\alpha$-stable random variables under
sublinear expectation. The foundation of our proof is an interior
regularity estimate for partial integro-differential equations (PIDEs).
A classical generalized central limit theorem is recovered as a special
case, provided a mild but natural additional condition holds. Our
approach contrasts with previous arguments for the result in the linear
setting which have typically relied upon tools
that are non-existent in the sublinear framework, for example,
characteristic functions.
\end{abstract}

%
\end{frontmatter}

\section{Introduction}\label{intro}

The purpose of this manuscript is to prove a generalized central limit
theorem for $\alpha$-stable random variables in the setting of
sublinear expectation. Such a result complements the limit theorems for
$G$-normal random variables due to Peng and others in this context and
answers in the affirmative a question posed by Neufeld and Nutz in
\cite{nutz+neufeld} (see below).

When working with a sublinear expectation, one is simultaneously
considering a potentially uncountably infinite and non-dominated
collection of probability measures. A construction of this kind is
motivated by the study of pricing under volatility uncertainty.
Needless to say, a variety of frequently called upon devices from the
classical setting are unavailable. The complications encompass further
issues as well: new behaviors are occasionally observed like those
outlined in \cite{bay+munk}.

Analogues of significant theorems from classical probability and
stochastic analysis are nevertheless moderately abundant. For instance,
versions of the law of large numbers can be found in \cite
{peng+lln+clt} and \cite{peng+clt}; the martingale representation
theorem is given in \cite{soner+touzi+zhang,song+mart} and
\cite{peng+song+zhang};  Girsanov's theorem is obtained
in \cite{xu+shang+zhang,osuka} and \cite
{hu+ji+peng+song}; and a Donsker-type result is shown in \cite
{dolinsky+nutz+soner}. To conduct investigations along these lines,
standard proofs must often be re-imagined. For instance, Peng's proof
of the central limit theorem under sublinear expectation in \cite
{peng+lln+clt} resorts to interior regularity estimates for fully
nonlinear parabolic partial differential equations (PDEs). His idea has
since been extended to prove a number of variants of his original
result, for example, see \cite{peng+clt,li+shi,hu+zhou+clt} and \cite{zhang+chen}.

We will operate in the sublinear expectation framework unless
explicitly indicated otherwise. The objects of our special attention
here, the $\alpha$-stable random variables for $\alpha\in ( 1 ,
2  )$, were introduced in \cite{nutz+neufeld}. The authors
pondered whether or not these could be the subject of a generalized
central limit theorem. Classical generalized central limit theorems
ordinarily come in one of three flavors:
\begin{enumerate}[(iii)]
\item[(i)] a statement indicating that a random variable has a non-empty
domain of attraction if and only if it is $\alpha$-stable such as
Theorem~2.1.1 in \cite{ibragimov+linnik},
\item[(ii)] a characterization theorem for the domain of attraction of an
$\alpha$-stable random variable such as Theorem~2.6.1 in \cite
{ibragimov+linnik}, or
\item[(iii)] a characterization theorem for the domain of normal attraction
for an $\alpha$-stable random variable such as Theorem~2.6.7 in \cite
{ibragimov+linnik}.
\end{enumerate}
Recall that an i.i.d. sequence $ (Y_i  )_{i=1}^{\infty}$ of
random variables is in the {\it domain of attraction} of a random
variable $X$ if there exist sequences of constants $ (A_i
)_{i=1}^{\infty}$ and $ (B_i  )_{i=1}^{\infty}$ so that
\[
B_n \sum_{i = 1}^n
Y_i - A_n
\]
converges in distribution to $X$ as $n \rightarrow\infty$. $ (Y_i  )_{i=1}^{\infty}$ is in the {\it domain of normal
attraction} of $X$ if
\[
B_n = \frac{1}{b n^{1/\alpha}}
\]
for some $b > 0$.

We confine our search to the direction suggested by (iii) because of
the particular importance classically of results of this type (cf. the
central limit theorem). Our main findings are summarized in Theorem~\ref{newmr}, which details sufficient conditions for membership in
the domain of normal attraction of a given $\alpha$-stable random
variable. While the initial appearance of our distributional hypotheses
is perhaps forbidding, in point of fact, our assumptions are
manageable. This is illustrated by the discussion immediately following
Theorem~\ref{newmr}, as well as Examples~\ref{x1examp} and~\ref
{bignlexmp}.

Example~\ref{x1examp} establishes that the $\alpha$-stable random
variables under consideration are in their own domain of normal
attraction. Although one need not apply Theorem~\ref{newmr} for this
purpose, the write-up serves a clarifying role and any credible result
clearly must pass this litmus test.

Example~\ref{bignlexmp} is more substantive. Setting aside a few
mild ``uniformity'' conditions which arise due to the supremum, this
example can be understood in an intuitive manner (see Section~\ref{sectexmp}). This falls out of our analysis just below Theorem~\ref
{newmr}, where we describe the relationship between our work and the
classical result noted in (iii) above. More specifically, Theorem~\ref
{newmr} detects all classical random variables in this collection with
mean zero and a cumulative distribution function (cdf) that satisfies a
small differentiability requirement. An extra regularity condition on
the cdf is unavoidable, as one must translate its form into properties
that can be stated only in terms of expectation.

The strategy of our proof is to reduce demonstrating convergence in
distribution to showing that a certain limit involving the solution to
the backward version of our generating PIDE is zero. Upon breaking up
our domain and summing the corresponding increments of the solution,
regularity properties of this function are employed to argue that size
of the terms being added together decay rapidly enough in the limit to
furnish the desired conclusion. This general scheme is similar to that
initiated in \cite{peng+lln+clt}, except that the generating equation
there is
\begin{eqnarray*}
\partial_t u - \frac{1}{2} \bigl( \overline{
\sigma}^2 ( \partial_{xx} u )^+ - \underline{
\sigma}^2 (\partial _{xx} u )^- \bigr)&=& 0 ,\qquad (t, x )
\in (0 , \infty ) \times\mathbb{R},
\\
u ( 0 ,x ) &=& \psi ( x ),\qquad x \in\mathbb{R}
\end{eqnarray*}
for some $0 \leq\underline{\sigma}^2 \leq\overline{\sigma}^2$ and
appropriate function $\psi$. Recall that this equation is known as the
Barenblatt equation if $\underline{\sigma}^2 > 0$ and has been
studied in \cite{baren} and \cite{avell+levy+paras}, for instance.
Ours is given by (\ref{uPIDEback}), a difference that leads to a few
difficulties as reflected by the increased complexity of our
hypotheses. To overcome these difficulties, we use the technology from
\cite{lara+dav+crit+dr,lara+dav+sig+al} and \cite{caff+cabre}.

The work in this paper offers a step toward understanding $\alpha
$-stability under sublinear expectation. The simple interpretation
admitted by Example~\ref{bignlexmp} is promising, as developing
intuition in this environment is usually a tough undertaking for the
reasons mentioned previously.

A brief overview of necessary background material can be found in
Section~\ref{sectback}. We prove our main result and discuss its
connection to the classical case in Section~\ref{mrsection}. Examples
highlighting the applications of our main result are contained in
Section~\ref{sectexmp}. We give some prerequisite material for the proof of the
essential interior regularity estimate for our PIDE in
Appendix~\ref{PIDEbackgroundapp}. The proof of this estimate is in Appendix~\ref{PIDEappend}.

\section{Background}\label{sectback}
We now offer a concise account of those aspects of sublinear
expectations, $\alpha$-stable random variables, and PIDEs which are
required for the sequel.\footnote{Further information on PIDE
interior regularity theory is contained in Appendix~\ref{PIDEbackgroundapp}.} References for more comprehensive treatments of
these topics are also included for the convenience of the interested reader.

\begin{defn}
Let $\mathcal{H}$ be a collection of real-valued functions on a set
$\Omega$. A {\it sublinear expectation} is an operator $\mathcal{E} :
\mathcal{H} \longrightarrow\mathbb{R}$ which is:
\begin{enumerate}[(iii)]
\item[(i)] monotonic: $\mathcal{E}  [X  ] \leq\mathcal{E}
 [Y  ]$ if $X \leq Y$,
\item[(ii)] constant-preserving: $\mathcal{E}  [c  ] = c$ for any
$c\in\mathbb{R}$,
\item[(iii)] sub-additive: $\mathcal{E}  [ X + Y ] \leq \mathcal
{E}  [X  ] + \mathcal{E}  [Y  ]$, and
\item[(iv)] positive homogeneous: $\mathcal{E}  [ \lambda X ] =
\lambda\mathcal{E}  [X  ]$ for $\lambda\geq0$.
\end{enumerate}
The triple $ ( \Omega, \mathcal{H} , \mathcal{E}  )$ is
called a {\it sublinear expectation space}.
\end{defn}

One views $\mathcal{H}$ as a space of random variables on $\Omega$.
Typically, it is assumed that $\mathcal{H}$
\begin{enumerate}[(iii)]
\item[(i)] is a linear space,
\item[(ii)] contains all constant functions, and
\item[(iii)] contains $\psi ( X_1 , X_2, \dots, X_n  )$ for every
$X_1 , X_2, \dots, X_n \in\mathcal{H}$ and $\psi\in C_{\mathrm{b.Lip}}  ( \mathbb{R}^n  )$, where $C_{\mathrm{b.Lip}}  ( \mathbb{R}^n
 )$ is the set of bounded Lipschitz functions on $\mathbb{R}^n$;
\end{enumerate}
however, we will expend little attention on either $\Omega$ or
$\mathcal{H}$. Delicacy needs to be exercised while computing
sublinear expectations. A rare instance when a classical technique can
be justly employed is the following.

\begin{lem}\label{lemsubexp}
Consider two random variables $X, Y \in\mathcal{H}$ such that
$\mathcal{E}  [Y  ] = - \mathcal{E}  [-Y  ]$. Then
\[
\mathcal{E} [X + \alpha Y ] = \mathcal{E} [X ] + \alpha\mathcal{E} [Y ]
\]
for all $\alpha\in\mathbb{R}$.
\end{lem}

This result is notably useful in the case where $\mathcal{E}  [Y
 ] = \mathcal{E}  [-Y  ] = 0$.

\begin{defn}
A random variable $Y \in\mathcal{H}$ is said to be {\it independent}
from a random variable $X \in\mathcal{H}$ if for all $\psi\in
C_{\mathrm{b.Lip}}  (\mathbb{R}^{2}  )$, we have
\[
\mathcal{E} \bigl[ \psi ( X , Y ) \bigr] = \mathcal{E} \bigl[ \mathcal{E} \bigl[
\psi ( x , Y ) \bigr]_{x =
X} \bigr] .
\]
\end{defn}

Observe the deliberate wording. This choice is crucial, as independence
can be asymmetric in our context. Note that this definition reduces to
the traditional one if $\mathcal{E}$ is a classical expectation. The
same is true for the next three concepts.

\begin{defn}
Let $X$, $Y$ and $ ( Y_n  )_{n=1}^{\infty}$ be random
variables, that is, $X$, $Y$ and $ ( Y_n  )_{n=1}^{\infty}
\in\mathcal{H}$.
\begin{enumerate}[(iii)]
\item[(i)]$X$ and $Y$ are {\it identically distributed}, denoted $X \sim
Y$, if
\[
\mathcal{E} \bigl[ \psi ( X ) \bigr] = \mathcal{E} \bigl[ \psi ( Y ) \bigr]
\]
for all $\psi\in C_{\mathrm{b.Lip}}  ( \mathbb{R}  )$.
\item[(ii)] If $X$ and $Y$ are identically distributed and $Y$ is independent
from $X$, then $Y$ is an {\it independent copy} of $X$.
\item[(iii)]$ ( Y_n  )_{n=1}^{\infty}$ {\it converges in
distribution} to $Y$, which we denote by $Y_n \xrightarrow{d} Y$, if
\[
\lim_{n \rightarrow\infty} \mathcal{E} \bigl[ \psi ( Y_n ) \bigr]
= \mathcal{E} \bigl[ \psi ( Y ) \bigr]
\]
for all $\psi\in C_{\mathrm{b.Lip}}  ( \mathbb{R}  )$.
\end{enumerate}
\end{defn}

Random variables need not be defined on the same space to have
appropriate notions of (i) or (iii). In this case, the above
definitions require the obvious notational modifications. Further
details concerning general sublinear expectation spaces can be found in
\cite{peng+ln} or \cite{peng+sc}.

\begin{defn}
Let $\alpha\in ( 0 , 2  ]$. A random variable $X$ is said
to be {\it(strictly) $\alpha$-stable} if for all $a$, $b \geq0$,
\[
a X + b Y
\]
and
\[
\bigl( a^{\alpha} + b^{\alpha} \bigr)^{1/ \alpha} X
\]
are identically distributed, where $Y$ is an independent copy of $X$.
\end{defn}

Three examples of $\alpha$-stable random variables exist in the
current literature. For $\alpha=1$, there are the maximal random
variables discussed in references such as \cite{peng+clt,peng+sc} and \cite{hu+li}. When $\alpha= 2$, we have the $G$-normal
random variables of Peng. Resources on this topic are plentiful and
include \cite{peng+ln,peng+survey,peng+sc,luo+jia} and \cite{bay+munk}. If $\alpha\in ( 1 , 2  )$,
we can consider $X_1$ for a nonlinear $\alpha$-stable L{\'e}vy process
$ ( X_t  )_{t \geq0}$ in the framework of \cite
{nutz+neufeld}. Our focus shall be restricted to the last situation.

The construction of nonlinear L{\'e}vy processes in \cite
{nutz+neufeld} extends that studied in \cite{peng+levy,ren,paczka+levy} and \cite{paczka+calc} and is much more general
than our present objectives demand. We limit our presentation to a few
key ideas. Let
\begin{enumerate}[(iii)]
\item[(i)]$\alpha\in (1, 2  )$;
\item[(ii)]$K_{\pm}$ be a bounded measurable subset of $\mathbb{R}_{+}$;
\item[(iii)]$F_{k_{\pm}}$ be the $\alpha$-stable L\'{e}vy measure
\[
F_{k_{\pm}} ( dz ) = ( k_{-} \mathbf{1}_{ (-
\infty, 0 )} +
k_{+} \mathbf{1}_{ (0, \infty )} ) (z )\llvert z\rrvert ^{-\alpha- 1}\,dz
\]
for all $k_{\pm} \in K_{\pm}$; and
\item[(iv)]$\Theta=  \{  ( 0 , 0, F_{k_{\pm}}  ) : k_{\pm
} \in K_{\pm}  \}$.
\end{enumerate}
One can then produce a process $ ( X_t  )_{t \geq0}$ which
is a nonlinear L{\'e}vy process whose local characteristics are
described by the set of L{\'e}vy triplets $\Theta$. This means the following:
\begin{enumerate}[(iii)]
\item[(i)]$ ( X_t  )_{t \geq0}$ is a real-valued c{\`a}dl{\`a}g process.
\item[(ii)]$X_0 = 0$.
\item[(iii)]$ ( X_t  )_{t \geq0}$ has stationary increments, that
is, $X_{t} - X_{s}$ and $X_{t-s}$ are identically distributed for all
$0 \leq s \leq t$.
\item[(iv)]$ ( X_t  )_{t \geq0}$ has independent increments,
that is, $X_{t} - X_{s}$ is independent from $ ( X_{s_1} ,
\dots, X_{s_n}  )$ for all $0 \leq s_1 \leq\cdots\leq s_n \leq
s \leq t$.
\item[(v)] If $\psi\in C_{\mathrm{b.Lip}}  ( \mathbb{R}  )$ and $u$ is
defined by
\[
u (t,x ) = \mathcal{E} \bigl[ \psi ( x + X_t ) \bigr]
\]
for all $ ( t, x ) \in[0 , \infty ) \times
\mathbb{R}$, then $u$ is the unique\footnote{The uniqueness
of a viscosity solution of (\ref{uPIDEback}) can be viewed as a
special case of Theorem~2.5 in \cite{nutz+neufeld}.} viscosity
solution\footnote{We take the following definition from
Section~2.2 of \cite{nutz+neufeld}. Let $C^{2,3}_b  (  ( 0 ,
\infty ) \times\mathbb{R}  )$ denote the set of
functions on \mbox{$ ( 0 , \infty ) \times\mathbb{R}$} having
bounded continuous partial derivatives up to the second and third order
in $t$ and $x$, respectively. A bounded upper semicontinuous function
$u$ on $ [ 0 , \infty ) \times\mathbb{R}$ is a {\it
viscosity subsolution} of (\ref{uPIDEback}) if
\[
u ( 0, \cdot ) \leq\psi ( \cdot )
\]
and for any $ ( t , x  ) \in ( 0 , \infty
 ) \times\mathbb{R}$,
\[
\partial_t \varphi ( t, x ) - \sup_{k_{\pm} \in K_{\pm}} \biggl
\{ \int_{\mathbb{R}} \delta_z \varphi ( t, x )
F_{k_{\pm}} (dz ) \biggr\} \leq0
\]
whenever $\varphi\in C^{2,3}_b  (  ( 0 , \infty
 ) \times\mathbb{R}  )$ is such that
\[
\varphi\geq u
\]
on \mbox{$ ( 0 , \infty ) \times\mathbb{R}$} and
\[
\varphi ( t, x ) = u ( t, x ).
\]
To define a {\it viscosity supersolution} of (\ref{uPIDEback}), one reverses the
inequalities and semicontinuity. A bounded
continuous function is a {\it viscosity solution} of (\ref{uPIDEback})
if it is both a viscosity subsolution and supersolution.
Viscosity solutions of the other PIDEs appearing in this paper, for
example, see Lemma~\ref{inithold}, are defined similarly.} of
\begin{eqnarray}
\label{uPIDEback} \partial_t u (t, x ) - \sup_{k_{\pm} \in
K_{\pm}}
\biggl\{ \int_{\mathbb{R}} \delta_z u ( t, x )
F_{k_{\pm}} (dz ) \biggr\}&=& 0 ,\qquad (t, x ) \in (0 , \infty ) \times
\mathbb{R}
\nonumber
\\[-8pt]
\\[-8pt]
u ( 0 ,x ) &=& \psi ( x ), \qquad x \in\mathbb{R}.
\nonumber
\end{eqnarray}
\end{enumerate}
Here, we use the notation
\[
\delta_z u ( t, x ) := u ( t , x + z ) - u ( t, x ) -
\partial_x u ( t , x ) z
\]
since the right-hand side of this equation as well as similar
expressions will frequently occur throughout the paper.

A critical feature of this setup is that if $\Theta$ is a singleton,
$ ( X_t  )_{t \geq0}$ is a classical L{\'e}vy process with
triplet $\Theta$. That $X_1$ actually is an $\alpha$-stable random
variable is not immediately obvious. We give a brief argument in
Example~\ref{x1examp}, but the core of this observation is a result
from \cite{nutz+neufeld} (see Example~2.7).

\begin{lem}\label{stabprocprop}
For all $\beta>0$ and $t \geq0$, $X_{\beta t}$ and $\beta^{1/ \alpha
} X_{t}$ are identically distributed.
\end{lem}

The dynamic programming principle in Lemma~\ref{nnlem51}
(see Lemma~5.1 in \cite{nutz+neufeld}) and the absolute value bound
in Lemma~\ref{absvalprop} (see Lemma~5.2 in \cite
{nutz+neufeld}) also play a central role when using our main result to
check that $X_1$ is in its own domain of normal attraction.

\begin{lem}\label{nnlem51}
For all $0 \leq s \leq t< \infty$ and $x \in\mathbb{R}$,
\[
u ( t , x ) = \mathcal{E} \bigl[ u ( t - s , x + X_s ) \bigr].
\]
\end{lem}

\begin{lem}\label{absvalprop}
We have that
\[
\mathcal{E} \bigl[\llvert X_{1} \rrvert \bigr] < \infty.
\]
\end{lem}

The remaining essential ingredients for our purposes describe the
regularity of $u$. 
The first result describes properties of $u$
which are valid on the whole domain. It is a special case of Lemma~5.3
in \cite{nutz+neufeld}.

\begin{lem}\label{nnlem53}
The function $u$ is uniformly bounded by $\llVert   \psi\rrVert
_{L^{\infty}  ( \mathbb{R}  )}$ and jointly continuous.
More precisely, $u  ( t , \cdot )$ is Lipschitz continuous
with constant $\operatorname{Lip} ( \psi )$, the Lipschitz constant of
$\psi$, and $u  ( \cdot, x  )$ is locally $1/2$-H{\"o}lder
continuous with a constant depending only on $\operatorname{Lip} (\psi )$ and
\[
\sup_{k_{\pm} \in K_{\pm}} \biggl\{ \int_{\mathbb
{R}} \llvert z
\rrvert \wedge\llvert z\rrvert ^2 F_{k_{\pm}} (dz ) \biggr\} <
\infty.
\]
\end{lem}

We will require
even stronger regularity estimates for $u$. To obtain these, we must
restrict our attention to the interior of the domain.

\begin{prop}\label{uPIDEregprop}
Suppose that for some $\lambda$, $\Lambda> 0$, we know $\lambda<
k_{\pm} < \Lambda$ for all $k_{\pm} \in K_{\pm}$. For any $h >0$,
\begin{enumerate}[(iii)]
\item[(i)]$\partial_t u$ and $\partial_x u$ exist and are bounded on
$ [h, h+1  ] \times\mathbb{R}$;
\item[(ii)] there are constants $C$, $\gamma>0$ such that
\begin{eqnarray*}
\bigl\llvert \partial_t u ( t_0 , x ) -
\partial_t u ( t_1 , x ) \bigr\rrvert &\leq& C \llvert
t_0 - t_1 \rrvert ^{\gamma/
\alpha},
\\
\bigl\llvert \partial_t u ( t , x_0 ) -
\partial_t u ( t , x_1 ) \bigr\rrvert &\leq& C \llvert
x_0 - x_1 \rrvert ^{\gamma}
\end{eqnarray*}
for all $ ( t_0 , x  )$, $ ( t_1 , x  )$, $ ( t , x_0  )$, $ ( t , x_1  ) \in [h, h+1
] \times\mathbb{R}$;
\item[(iii)]$u$ is a classical solution of (\ref{uPIDEback}) on $ [h,
h+1  ] \times\mathbb{R}$; and
\item[(iv)] if $K_{\pm}$ contains exactly one pair $ \{ k_{\pm}
\}$, then $\partial_{xx}^2 u$ exists and is bounded on $ [h, h+1
 ] \times\mathbb{R}$.
\end{enumerate}
\end{prop}

The proof of Proposition~\ref{uPIDEregprop} can be found in
Appendix~\ref{PIDEappend}.

\section{Main result}\label{mrsection}
To facilitate our discussion in the sequel, we now fix some notation.
Compared with Section~\ref{sectback}, we make only one alteration to
our nonlinear $\alpha$-stable L{\'e}vy process $ ( X_t
)_{t \geq0}$: additionally assume that $K_{\pm}$ is a subset of
$ ( \lambda, \Lambda )$ for some $\lambda$, $\Lambda>
0$. We will make use of this in conjunction with Proposition~\ref{uPIDEregprop}.

We also consider a sequence $ (Y_i  )_{i = 1}^{\infty}$ of
random variables on some sublinear expectation space. The only aspect
of this space that we will invoke directly is the sublinear expectation
itself, say~$\mathcal{E}^{\prime}$. Distinguishing between $\mathcal
{E}$ and $\mathcal{E}^{\prime}$ will be convenient for Example~\ref
{bignlexmp}. We further specify that $ (Y_i  )_{i =
1}^{\infty}$ is i.i.d. in the sense that $Y_{i+1}$ is independent from
$ ( Y_1 , Y_2, \ldots, Y_i  )$ and $Y_{i+1} \sim Y_i$ for
all $i \geq1$. After proper normalization,
\[
S_n := \sum_{i = 1}^n
Y_i
\]
will be the sequence attracted to $X_1$.

\begin{thmo}\label{newmr}
Suppose that
\begin{enumerate}[(iii)]
\item[(i)]$\mathcal{E}^{\prime}  [ Y_1  ] = \mathcal
{E}^{\prime}  [ - Y_1  ] = 0$;
\item[(ii)]$\mathcal{E}^{\prime}  [ \llvert   Y_1 \rrvert    ] <
\infty$; and
\item[(iii)] for any $0 < h < 1$ and $\psi\in C_{\mathrm{b.Lip}}  ( \mathbb{R}
 )$,
\begin{equation}
\label{disthyp} n \biggl\llvert \mathcal{E}^{\prime} \bigl[
\delta_{B_n Y_1} v ( t, x ) \bigr] - \biggl( \frac{1}{n} \biggr) \sup
_{k_{\pm} \in K_{\pm}} \biggl\{ \int_{\mathbb{R}}
\delta_z v ( t, x ) F_{k_{\pm}} (dz ) \biggr\} \biggr\rrvert
\rightarrow0
\end{equation}
uniformly on $ [ 0 , 1 ] \times\mathbb{R}$ as $n
\rightarrow\infty$, where $v$ is the unique viscosity solution of
\begin{eqnarray}
\label{vPIDE} \partial_t v (t, x ) + \sup_{k_{\pm} \in
K_{\pm}}
\biggl\{ \int_{\mathbb{R}} \delta_z v ( t, x )
F_{k_{\pm}} (dz ) \biggr\} &=& 0 ,\qquad (t, x ) \in (-h , 1+h ) \times
\mathbb{R}
\nonumber
\\[-8pt]
\\[-8pt]
v ( 1+ h, x ) &=& \psi ( x ) ,\qquad x \in \mathbb{R}.
\nonumber
\end{eqnarray}
\end{enumerate}
Then
\[
B_n S_n \xrightarrow{d} X_1
\]
as $n \rightarrow\infty$.
\end{thmo}

Admittedly, a cursory glance over our hypotheses leaves one with the
impression that they are intractable. The opposite is true. Before
presenting the proof of Theorem~\ref{newmr}, let us demonstrate that
when our attention is confined to the classical case, we are imposing
only a mild and natural supplementary restriction on the attracted
random variable. In addition to being a significant remark in itself,
this work also underlies Example~\ref{bignlexmp}.

Assume that $\Theta$ is a singleton. Since $ ( X_t  )_{t
\geq0}$ is the classical L\'{e}vy process with triplet $ ( 0 , 0
, F_{k_{\pm}}  )$, the characteristic function of $X_1$, denoted
$\varphi_{X_1}$, is given by
\[
\varphi_{X_1} ( t ) = \exp \biggl( k_{-} \int
_{- \infty}^0 \frac{ \exp ( itz  ) - 1 - itz }{\llvert   z \rrvert  ^{\alpha+1}}\,dz + k_{+}
\int_{0}^{\infty} \frac{
\exp ( itz  ) - 1 - itz }{ z^{\alpha+1}}\,dz \biggr)
\]
for all $t \in\mathbb{R}$. In the case where $Y_1$ is a classical
random variable with mean zero, Theorem~2.6.7 from \cite
{ibragimov+linnik} implies that
\[
B_n S_n \xrightarrow{d} X_1
\]
as $n \rightarrow\infty$ if and only if the cdf of $Y_1$, denoted
$F_{Y_1}$, has the form
\[
F_{Y_1} ( z ) = \cases{ \bigl[ b^{\alpha} ( k_- / \alpha ) +
\beta_1 ( z ) \bigr] \dfrac{1}{\llvert  z \rrvert  ^\alpha}, & $z < 0$,
\cr
1 - \bigl[
b^{\alpha} ( k_+ / \alpha ) + \beta_2 ( z ) \bigr]
\dfrac{1}{z^\alpha}, & $z > 0$, } %
\]
for some functions $\beta_1$ and $\beta_2$ satisfying
\[
\lim_{z \rightarrow- \infty} \beta_1 ( z ) = \lim
_{z
\rightarrow\infty} \beta_2 ( z ) = 0 .
\]

As there is no appropriate counterpart of the cdf in the sublinear
setting, we must recast this condition using expectation. To do so
requires $F_{Y_1}$ to possess further regularity properties. For
convenience, say that after an extension, the $\beta_i$'s are
continuously differentiable on their respective closed half-lines. This
is the lone extra requirement we shall need.

It follows that
\[
\mathbb{E} \bigl[ \llvert Y_1 \rrvert \bigr] < \infty
\]
since
\begin{eqnarray}
\label{absvalint1} \int_{0}^{\infty} z \,
dF_{Y_1} ( z ) &=& - \int_{0}^{1}
\frac{\beta_2^{\prime}  ( z  )}{z^{\alpha-1}}\,dz + \int_{0}^{1}
\frac{b^{\alpha} k_+ + \alpha\beta_2  ( z
)}{z^{\alpha}}\,dz + \beta_2 ( 1 )
\nonumber
\\
&&{} + \int_{1}^{\infty} \frac{\beta_2  ( z
)}{z^{\alpha}}\,dz + \int
_{1}^{\infty} \frac{b^{\alpha}
k_+}{z^{\alpha}}\,dz
\\
&<& \infty\nonumber
\end{eqnarray}
and similarly for the integral along the negative half-line. One could
have cited Theorem~2.6.4 of \cite{ibragimov+linnik} instead, but (\ref
{absvalint1}) will be helpful in Example~\ref{bignlexmp}. We also get
\begin{eqnarray}
\label{linequivdisthyp2}&& n \biggl\llvert \mathbb{E} \bigl[ \delta_{B_n Y_1} v ( t,
x ) \bigr] - \biggl( \frac{1}{n} \biggr) \int_{\mathbb{R}}
\delta _z v ( t, x ) F_{k_{\pm}} (dz ) \biggr\rrvert
\nonumber
\\
&&\quad = \biggl( \frac{1}{b^\alpha} \biggr) \bigg\vert \int_{\mathbb{R}}
\delta_z v ( t, x ) \biggl( \frac{ \beta
_1^{\prime}  ( B_{n}^{-1} z  ) \llvert  B_{n}^{-1} z \rrvert
+ \alpha\beta_1  ( B_{n}^{-1} z  ) }{ \llvert   z \rrvert  ^{\alpha+ 1}} {
\bf1}_{ ( -\infty, 0  )} ( z )
\\
&&\qquad {} + \frac{- \beta_2^{\prime}
 ( B_{n}^{-1} z  ) \llvert  B_{n}^{-1} z \rrvert   + \alpha
\beta_2  ( B_{n}^{-1} z  ) }{ \llvert   z \rrvert  ^{\alpha+
1}} \mathbf{1}_{ ( 0 , \infty )} ( z ) \biggr)\,dz \bigg\vert
\nonumber
\end{eqnarray}
for all $ ( t , x  ) \in [0, 1 ] \times\mathbb
{R}$ and $n \geq1$ by changing variables.

A careful application of elementary estimates shows that this last
expression tends to zero uniformly on $ [0, 1 ] \times
\mathbb{R}$ as $n \rightarrow\infty$. To see this, note that we can
choose an upper bound, say~$M_1$, for $\llvert   \partial_{xx} v \rrvert  $, $\llvert   \partial_{x} v \rrvert  $, and $\llvert   v \rrvert  $ on
$ [0, 1 ] \times\mathbb{R}$ by Lemma~\ref{nnlem53} and
Proposition~\ref{uPIDEregprop}. Then using integration by parts and
the dominated convergence theorem,
\begin{eqnarray}
\label{intbd1} && \biggl\llvert \int_{1}^{\infty}
\delta_z v ( t, x ) \biggl( \frac{- \beta_2^{\prime}  ( B_{n}^{-1} z  ) \llvert  B_{n}^{-1} z \rrvert   + \alpha\beta_2  ( B_{n}^{-1} z  )
}{ \llvert   z \rrvert  ^{\alpha+ 1}} \biggr)\,dz \biggr
\rrvert
\nonumber
\\
&&\quad = \biggl\llvert \delta_{1} v ( t, x ) \beta_2
\bigl( B_{n}^{-1} \bigr)\nonumber
\\
&&\qquad {}+ \int_{1}^{\infty}
\frac{\beta_2  (
B_{n}^{-1} z  ) }{z^{\alpha}} \bigl[ \partial_xv (t , x + z ) -
\partial_x v (t, x ) \bigr]\,dz \biggr\rrvert
\\
&&\quad \leq3 M_1 \bigl\llvert \beta_2 \bigl(
B_{n}^{-1} \bigr) \bigr\rrvert + 2 M_1 \int
_{1}^{\infty} \frac{\llvert  \beta_2  ( B_{n}^{-1} z
 ) \rrvert  }{z^{\alpha}}\,dz
\nonumber
\\
&&\quad \rightarrow0
\nonumber
\end{eqnarray}
as $n \rightarrow\infty$. The mean value theorem and a change of
variables give
\begin{eqnarray}
\label{intbd2} && \biggl\llvert \int_{0}^{B_n}
\delta_z v ( t, x ) \biggl( \frac{- \beta_2^{\prime}  ( B_{n}^{-1} z  ) \llvert  B_{n}^{-1} z \rrvert   + \alpha\beta_2  ( B_{n}^{-1} z  )
}{ \llvert   z \rrvert  ^{\alpha+ 1}} \biggr)\,dz \biggr
\rrvert
\nonumber
\\
&&\quad \leq \int_{0}^{B_n} M_1
\frac{\llvert  {-} \beta_2^{\prime}  (
B_{n}^{-1} z  )  ( B_{n}^{-1} z  ) + \alpha\beta_2
 ( B_{n}^{-1} z  ) \rrvert   }{ z^{\alpha- 1}}\,dz
\\
&&\quad = \biggl( \frac{ M_1}{ b^{2-\alpha} n^{ \sfrac{2}{\alpha} -1 }} \biggr) \int_{0}^{1}
\frac{\llvert   {-} \beta_2^{\prime}  ( z
 ) z + \alpha\beta_2  ( z  )\rrvert   }{ z^{\alpha-
1}}\,dz
\nonumber
\\
&&\quad \rightarrow0
\nonumber
\end{eqnarray}
as $n \rightarrow\infty$. We have
\begin{eqnarray}
\label{intbd3} && \biggl\llvert \int_{B_n}^{1}
\delta_z v ( t, x ) \biggl( \frac{\alpha\beta_2  ( B_{n}^{-1} z  ) }{ \llvert   z
\rrvert  ^{\alpha+ 1}} \biggr)\,dz \biggr
\rrvert
\nonumber
\\
&&\quad \leq \int_{B_n}^{1} M_1
\frac{\llvert   \alpha\beta_2  (
B_{n}^{-1} z  )\rrvert   }{ z^{\alpha- 1}}\,dz
\nonumber
\\[-8pt]
\\[-8pt]
&&\quad \leq M_1 \alpha\int_{0}^{1}
\frac{\llvert   \beta_2  (
B_{n}^{-1} z  )\rrvert   }{ z^{\alpha- 1}}\,dz
\nonumber
\\
&&\quad \rightarrow0
\nonumber
\end{eqnarray}
as $n \rightarrow\infty$ by the mean value theorem and dominated
convergence theorem. Finally,
\begin{eqnarray}
\label{intbd4} && \biggl\llvert \int_{B_n}^{1}
\delta_z v ( t, x ) \biggl( \frac{- \beta_2^{\prime}  ( B_{n}^{-1} z  )  (B_{n}^{-1} z  ) }{ \llvert   z \rrvert  ^{\alpha+ 1}} \biggr)\,dz \biggr
\rrvert
\nonumber
\\
&&\quad = \bigg\vert {-} \delta_{1} v ( t, x ) \beta_2
\bigl(B_{n}^{-1} \bigr) + \delta_{B_n} v ( t, x ) (
B_n )^{-\alpha} \beta_2 ( 1 )
\nonumber
\\
&&\qquad {} + \int_{B_n}^{1} \bigl[
\partial_x v (t , x + z ) - \partial_x v (t, x ) \bigr]
\biggl( \frac{ \beta_2
 ( B_{n}^{-1} z  ) }{ z^{\alpha}} \biggr)\,dz
\nonumber
\\
&&\qquad {} - \alpha\int_{B_n}^{1}
\delta_z v ( t, x ) \biggl( \frac{ \beta_2  ( B_{n}^{-1} z  ) }{ z^{\alpha
+1}} \biggr)\,dz \bigg\vert
\nonumber
\\[-8pt]
\\[-8pt]
&&\quad \leq3 M_1 \bigl\llvert \beta_2
\bigl(B_{n}^{-1} \bigr) \bigr\rrvert + M_1 \bigl
\llvert \beta_2 (1 ) \bigr\rrvert \biggl( \frac
{1}{b^{2-\alpha} n^{\sfrac{2}{\alpha} - 1}} \biggr) +
\int_{B_n}^{1} M_1 \frac{ \llvert   \beta_2  ( B_{n}^{-1} z  )
\rrvert  }{ z^{\alpha-1}}
\,dz
\nonumber
\\
&&\qquad {} + \alpha\int_{B_n}^{1}
M_1 \frac{ \llvert   \beta_2  (
B_{n}^{-1} z  ) \rrvert   }{ z^{\alpha-1}}\,dz
\nonumber
\\
&&\quad \leq3 M_1 \bigl\llvert \beta_2
\bigl(B_{n}^{-1} \bigr) \bigr\rrvert + M_1 \bigl
\llvert \beta_2 (1 ) \bigr\rrvert \biggl( \frac{1}{b^{2-\alpha
} n^{\sfrac{2}{\alpha} - 1}} \biggr) +
2 \alpha M_1 \int_{0}^{1}
\frac{\llvert   \beta_2  ( B_{n}^{-1} z  )\rrvert   }{
z^{\alpha- 1}}\,dz
\nonumber
\\
&&\quad \rightarrow0
\nonumber
\end{eqnarray}
as $n \rightarrow\infty$ by integration by parts, the dominated
convergence theorem, and the mean value theorem. The integrals along
the negative half-line are handled similarly.

Having established the connection between Theorem~\ref{newmr} and the
classical case, we finally present its proof.

\begin{pf*}{Proof of Theorem~\ref{newmr}}
We need to show that
\begin{equation}
\label{convdist} \lim_{n \rightarrow\infty} \mathcal{E}^{\prime} \bigl[
\psi ( B_n S_n ) \bigr] = \mathcal{E} \bigl[ \psi (
X_1 ) \bigr]
\end{equation}
for all $\psi\in C_{\mathrm{b.Lip}}  ( \mathbb{R}  )$. Our initial
step will be to reduce proving (\ref{convdist}) to proving (\ref
{hardlemlimit}). The advantage of doing so is that we can then
incorporate the regularity properties described in Lemma~\ref{nnlem53} and Proposition~\ref{uPIDEregprop}. These properties alone do
much of the heavy lifting in the estimates at the heart of the
argument, and our distributional assumptions do the rest.

Let $\psi\in C_{\mathrm{b.Lip}}  ( \mathbb{R}  )$, and define $u$ by
\begin{equation}
\label{defu} u (t,x ) = \mathcal{E} \bigl[ \psi ( x + X_t )
\bigr]
\end{equation}
for all $ ( t, x ) \in [0 , \infty ) \times
\mathbb{R}$. We know from Section~\ref{sectback} that $u$ is the
unique viscosity solution of (\ref{uPIDEback}).

It will be more convenient for our purposes to work with the backward
equation. Since we will soon rely on the interior regularity results of
Proposition~\ref{uPIDEregprop}, we also let $0 < h < 1$ and define
$v$ by
\begin{equation}
\label{relatvu} v (t,x ) = u (1 + h - t , x )
\end{equation}
for $ (t,x ) \in( -h , 1 + h ] \times\mathbb
{R}$. Then $v$ will be the unique viscosity solution of (\ref{vPIDE}).

Observe that $v$ inherits key regularity properties from $u$. At the
moment, it is enough to note that for any $ (t,x ) \in ( -h , 1 + h ] \times\mathbb{R}$, $v  ( \cdot, x )$
is $1/2$-H{\"o}lder continuous with some constant $K_1$ and $v  ( t
, \cdot )$ is Lipschitz continuous with constant $\operatorname{Lip} (
\psi )$ by Lemma~\ref{nnlem53}. Because the $t$-domain has
length $1+2h$ and $0 < h < 1$, the $1/2$-H{\"o}lder continuity is
uniform, and we can assume that $K_1$ does not depend on $h$. It
follows by (\ref{defu}) and (\ref{relatvu}) that
\begin{eqnarray*}
& &\limsup_{n \rightarrow\infty} \bigl\llvert \mathcal {E}^{\prime}
\bigl[ \psi ( B_n S_n ) \bigr] - \mathcal {E} \bigl[ \psi
( X_1 ) \bigr] \bigr\rrvert
\\
&&\quad \leq \limsup_{n \rightarrow\infty} \bigl( \bigl\llvert
\mathcal{E}^{\prime} \bigl[ \psi ( B_n S_n ) \bigr]
- \mathcal{E}^{\prime} \bigl[ v ( 1 , B_n S_n )
\bigr] \bigr\rrvert + \bigl\llvert \mathcal{E}^{\prime} \bigl[ v ( 1 ,
B_n S_n ) \bigr] - v (0, 0 ) \bigr\rrvert
\\
&&\qquad {} + \bigl\llvert v (0, 0 ) - \mathcal{E} \bigl[ \psi ( X_1
) \bigr] \bigr\rrvert \bigr)
\\
&&\quad = \limsup_{n \rightarrow\infty} \bigl( \bigl\llvert
\mathcal{E}^{\prime} \bigl[ v ( 1 + h , B_n S_n )
\bigr] - \mathcal{E}^{\prime} \bigl[ v ( 1 , B_n S_n
) \bigr] \bigr\rrvert + \bigl\llvert \mathcal{E}^{\prime} \bigl[ v ( 1 ,
B_n S_n ) \bigr] - v (0, 0 ) \bigr\rrvert
\\
&&\qquad {} + \bigl\llvert v (0, 0 ) - v (h , 0 ) \bigr\rrvert \bigr)
\\
&&\quad \leq \limsup_{n \rightarrow\infty} \bigl( \mathcal {E}^{\prime} [
K_1 \sqrt{h} ] + \bigl\llvert \mathcal {E}^{\prime} \bigl[ v (
1 , B_n S_n ) \bigr] - v (0, 0 ) \bigr\rrvert \bigr)+
K_1 \sqrt{h}
\\
&&\quad = 2 K_1 \sqrt{h} + \limsup_{n \rightarrow\infty} \bigl\llvert
\mathcal{E}^{\prime} \bigl[ v ( 1 , B_n S_n ) \bigr]
- v (0, 0 ) \bigr\rrvert .
\end{eqnarray*}
As $h$ is arbitrary, it is sufficient to show that
\begin{equation}
\label{hardlemlimit} \lim_{n \rightarrow\infty} \mathcal{E}^{\prime} \bigl[ v
( 1 , B_n S_n ) \bigr] = v (0, 0 ) .
\end{equation}
The required estimates are intricate, so we will give them in Lemma~\ref{hardlem} below.
\end{pf*}

\begin{lem}\label{hardlem}
In the setup of Theorem~\ref{newmr},
\[
\lim_{n \rightarrow\infty} \mathcal{E}^{\prime} \bigl[ v ( 1 ,
B_n S_n ) \bigr] = v (0, 0 ).
\]
\end{lem}

\begin{pf}
For all $n \geq3$,
\begin{eqnarray}
\label{initdecomp} &&v ( 1 , B_n S_n ) - v (0 , 0 )
\nonumber
\\
&&\quad = v ( 1 , B_n S_n ) - v \biggl(
\frac{n-1}{n} , B_n S_{n} \biggr) + \sum
_{i = 2}^{n-1} \biggl[ v \biggl( \frac{i}{n} ,
B_n S_{i+1} \biggr) - v \biggl( \frac{i-1}{n} ,
B_n S_{i} \biggr) \biggr]\quad
\\
&&\qquad {} + v \biggl( \frac{1}{n} , B_n S_{2}
\biggr) - v (0 , 0 ) .
\nonumber
\end{eqnarray}
Our analysis now becomes delicate. We would like to show that when we
apply $\mathcal{E}^\prime$ to (\ref{initdecomp}) and let $n
\rightarrow\infty$, the result goes to zero. Since the number of
terms in this decomposition is growing with $n$, we must prove that our
$v$-increments are decaying quite rapidly. The properties of $v$
arising from Lemma~\ref{nnlem53} are only enough to manage the first
and last terms. By the $1/2$-H{\"o}lder continuity of $v  ( \cdot,
x  )$,
\begin{equation}
\label{bd1} \mathcal{E}^{\prime} \biggl[ \biggl\llvert v ( 1 ,
B_n S_n ) - v \biggl( \frac{n-1}{n} ,
B_n S_{n} \biggr) \biggr\rrvert \biggr] \leq
\mathcal{E}^{\prime} \biggl[K_1 \sqrt{\frac{1}{n}} \biggr]
=K_1 \sqrt{\frac{1}{n}} .
\end{equation}
If we also use the Lipschitz continuity of $v  ( t, \cdot
)$ and the fact that $Y_2$ is independent from $Y_1$, we get
\begin{eqnarray}
\label{bd2} &&\mathcal{E}^{\prime} \biggl[ \biggl\llvert v \biggl(
\frac{1}{n} , B_n S_{2} \biggr) - v (0 , 0 ) \biggr
\rrvert \biggr]
\nonumber
\\
&&\quad \leq\mathcal{E}^{\prime} \biggl[ \biggl\llvert v \biggl(
\frac{1}{n} , B_n S_{2} \biggr) - v ( 0 ,
B_n S_{2} ) \biggr\rrvert \biggr] + \mathcal{E}^{\prime}
\bigl[ \bigl\llvert v ( 0 , B_n S_{2} ) - v (0 , 0 ) \bigr
\rrvert \bigr]
\nonumber
\\[-8pt]
\\[-8pt]
&&\quad \leq\mathcal{E}^{\prime} \biggl[K_1 \sqrt{
\frac{1}{n}} \biggr] + \mathcal{E}^{\prime} \bigl[ \operatorname{Lip}
(\psi ) B_n \llvert S_2 \rrvert \bigr]
\nonumber
\\
&&\quad \leq K_1 \sqrt{\frac{1}{n}} + 2 \operatorname{Lip} (
\psi ) B_n \mathcal{E}^{\prime} \bigl[ \llvert Y_1
\rrvert \bigr] .
\nonumber
\end{eqnarray}
We remark that although we only referred to $C_{\mathrm{b.Lip}}  ( \mathbb
{R}  )$ in our definition of independence, our manipulations are
still valid by Exercise 3.20 in \cite{peng+sc}.

Proposition~\ref{uPIDEregprop} allows us to control the remaining
terms. Again, this motivates our requirement that $K_{\pm} \subset
 ( \lambda, \Lambda )$ for some $0< \lambda< \Lambda$.
We can find a constant $K_2 > 0$ such that
$\partial_t v$ exists on $ [ 0, 1  ] \times\mathbb{R}$ and
\begin{eqnarray}
\label{reg} \bigl\llvert \partial_t v ( t_0 , x ) -
\partial_t v ( t_1 , x ) \bigr\rrvert &\leq&
K_2 \llvert t_0 - t_1 \rrvert ^{\gamma
/ \alpha},
\nonumber
\\[-8pt]
\\[-8pt]
\bigl\llvert \partial_t v ( t , x_0 ) -
\partial_t v ( t , x_1 ) \bigr\rrvert &\leq&
K_2 \llvert x_0 - x_1 \rrvert ^{\gamma}
\nonumber
\end{eqnarray}
for all $ ( t_0 , x  )$, $ ( t_1 , x  )$, $ ( t , x_0  )$, and $ ( t , x_1  ) \in [ 0, 1
 ] \times\mathbb{R}$.
We then break down the rest of (\ref{initdecomp}) a bit further. If
$2 \leq i \leq n-1$,
\begin{eqnarray*}
&&v \biggl( \frac{i}{n} , B_n S_{i+1} \biggr) - v
\biggl( \frac
{i-1}{n} , B_n S_{i} \biggr)
\\
&&\quad = v \biggl( \frac{i}{n} , B_n S_{i+1}
\biggr) - v \biggl( \frac
{i-1}{n} , B_n S_{i+1} \biggr)
- \partial_t v \biggl( \frac{i-1}{n} , B_n
S_{i} \biggr)\frac{1}{n}
\\
&&\qquad {} + \partial_t v \biggl( \frac{i-1}{n} ,
B_n S_{i} \biggr)\frac
{1}{n} + v \biggl(
\frac{i-1}{n} , B_n S_{i+1} \biggr) - v \biggl(
\frac{i-1}{n} , B_n S_i \biggr) .
\end{eqnarray*}
Let
\[
C^{n}_{i} = v \biggl( \frac{i}{n} ,
B_n S_{i+1} \biggr) - v \biggl( \frac{i-1}{n} ,
B_n S_{i+1} \biggr) - \partial_t v \biggl(
\frac
{i-1}{n} , B_n S_{i} \biggr)\frac{1}{n}
\]
and
\[
D^{n}_{i} = \partial_t v \biggl(
\frac{i-1}{n} , B_n S_{i} \biggr)\frac{1}{n} +
v \biggl( \frac{i-1}{n} , B_n S_{i+1} \biggr) - v
\biggl( \frac{i-1}{n} , B_n S_i \biggr) .
\]

We can establish an appropriate bound for the $C^{n}_{i}$'s using (\ref{reg}):
\begin{eqnarray*}
\bigl\llvert C^{n}_{i} \bigr\rrvert &=& \biggl\llvert
\frac{1}{n} \int_{0}^{1} \biggl[
\partial_t v \biggl( \frac{ i - 1+ \beta}{n} , B_n
S_{i + 1} \biggr) - \partial_t v \biggl( \frac{i-1}{n}
, B_n S_{i + 1} \biggr) \biggr]\,d\beta
\\
&&{} + \frac{1}{n} \biggl[ \partial_t v \biggl(
\frac
{i-1}{n} , B_n S_{i + 1} \biggr) -
\partial_t v \biggl( \frac
{i-1}{n} , B_n
S_{i} \biggr) \biggr] \biggr\rrvert
\\
&\leq&\frac{1}{n} \int_{0}^{1} \biggl
\llvert \partial_t v \biggl( \frac{
i - 1+ \beta}{n} , B_n
S_{i + 1} \biggr) - \partial_t v \biggl( \frac{i-1}{n}
, B_n S_{i + 1} \biggr) \biggr\rrvert \,d\beta
\\
&&{} + \frac{1}{n} \biggl\llvert \partial_t v \biggl(
\frac{i-1}{n} , B_n S_{i + 1} \biggr) -
\partial_t v \biggl( \frac{i-1}{n} , B_n
S_{i} \biggr) \biggr\rrvert
\\
&\leq& \frac{1}{n} \int_{0}^{1}
K_2 \biggl\llvert \frac{\beta}{n} \biggr\rrvert ^{\gamma/ \alpha}\,d
\beta+ \frac{1}{n} K_2 B_n^{\gamma} \llvert
Y_{i + 1} \rrvert ^{\gamma}
\\
&\leq&\frac{K_2 }{n} \biggl[ \biggl(\frac{1}{n} \biggr)^{\gamma/
\alpha}
+ B_n^{\gamma} \llvert Y_{i + 1} \rrvert
^{\gamma} \biggr] .
\end{eqnarray*}
Hence, for $2 \leq i \leq n-1$,
\begin{equation}
\label{bd3} \mathcal{E}^{\prime} \bigl[ \bigl\llvert
C^{n}_{i} \bigr\rrvert \bigr] \leq \frac{K_2 }{n}
\biggl[ \biggl(\frac{1}{n} \biggr)^{\gamma/ \alpha} + B_n^{\gamma}
\mathcal{E}^{\prime} \bigl[\llvert Y_{1} \rrvert ^{\gamma}
\bigr] \biggr]
\end{equation}
since $Y_{i+1}$ and $Y_1$ are identically distributed. Note that
hypothesis (ii) gives that
\[
\mathcal{E}^{\prime} \bigl[\llvert Y_{1} \rrvert
^{\gamma} \bigr] < \infty.
\]

While we need (\ref{reg}) to bound the $D^{n}_{i}$'s, we finally use
(\ref{disthyp}), too. Let $\varepsilon>0$. By (\ref{disthyp}), we can
find $N$ such that $n \geq N$ implies
\[
n \biggl\llvert \mathcal{E}^{\prime} \bigl[ \delta_{B_n Y_1} v ( t,
x ) \bigr] - \biggl( \frac{1}{n} \biggr) \sup_{k_{\pm} \in K_{\pm}}
\biggl\{ \int_{\mathbb{R}} \delta_z v ( t, x )
F_{k_{\pm}} (dz ) \biggr\} \biggr\rrvert < \varepsilon
\]
on $ [ 0 , 1 ] \times\mathbb{R}$. Now
\begin{eqnarray*}
&&\mathcal{E}^{\prime} \biggl[ v \biggl( \frac{i-1}{n} ,
B_n x + B_n Y_{1} \biggr) \biggr] - v \biggl(
\frac{i-1}{n} , B_n x \biggr)
\\
&&\quad  = \mathcal{E}^{\prime}
\biggl[ \delta_{B_n Y_1} v \biggl( \frac
{i-1}{n} , B_n x
\biggr) \biggr]
\end{eqnarray*}
by (i), so for these $n$,
\begin{eqnarray*}
&&n \biggl\llvert v \biggl( \frac{i-2}{n} , B_n x \biggr) -
\mathcal {E}^{\prime} \biggl[ v \biggl( \frac{i-1}{n} , B_n
x + B_n Y_{1} \biggr) \biggr] \biggr\rrvert
\\
&&\quad = n \bigg\vert v \biggl( \frac{i-2}{n} , B_n x \biggr) -
\mathcal {E}^{\prime} \biggl[ v \biggl( \frac{i-1}{n} , B_n
x + B_n Y_{1} \biggr) \biggr]
\\
& &\qquad {}+ v \biggl( \frac{i-1}{n} , B_n x \biggr) - v
\biggl( \frac
{i-1}{n} , B_n x \biggr)
\\
&&\qquad {} + \biggl( \frac{1}{n} \biggr) \partial_t v
\biggl( \frac
{i-1}{n} , B_n x \biggr) + \biggl( \frac{1}{n}
\biggr) \sup_{k_{\pm} \in K_{\pm}} \biggl\{ \int_{\mathbb{R}}
\delta_{z} v \biggl( \frac{i-1}{n} , B_n x \biggr)
F_{k_{\pm}} (dz ) \biggr\} \bigg\vert
\\
&&\quad \leq \biggl\llvert - \frac{ v  ( \vfrac{i-2}{n} , B_n x  ) - v
 ( \vfrac{i-1}{n} , B_n x  ) }{-1/n} + \partial_t v
\biggl( \frac{i-1}{n} , B_n x \biggr) \biggr\rrvert
\\
&&\qquad {} + n \biggl\llvert \mathcal{E}^{\prime} \biggl[
\delta_{B_n Y_1} v \biggl( \frac{i-1}{n} , B_n x \biggr)
\biggr]
\\
&&\qquad {}- \biggl( \frac{1}{n} \biggr) \sup_{k_{\pm} \in K_{\pm}} \biggl\{
\int_{\mathbb{R}} \delta_{z} v \biggl(
\frac{i-1}{n} , B_n x \biggr) F_{k_{\pm}} (dz ) \biggr\}
\biggr\rrvert
\\
&&\quad < \frac{K_2 }{n^{\gamma/ \alpha}} + \varepsilon
\end{eqnarray*}
by the mean value theorem, (\ref{vPIDE}), and (\ref{reg}). Then
\begin{eqnarray}
\label{bd4a} && \biggl\llvert \partial_t v \biggl(
\frac{i-1}{n} , B_n x \biggr)\frac
{1}{n} +
\mathcal{E}^{\prime} \biggl[ v \biggl( \frac{i-1}{n} , B_n
x + B_n Y_{i+1} \biggr) \biggr] - v \biggl(
\frac{i-1}{n} , B_n x \biggr) \biggr\rrvert
\nonumber
\\
&&\quad \leq\frac{1}{n} \biggl\llvert \partial_t v \biggl(
\frac{i-1}{n}, B_n x \biggr) + \frac{ v  ( \vfrac{i-2}{n} , B_n x  )- v  (
\vfrac{i-1}{n} , B_n x  )}{1/n} \biggr\rrvert
\nonumber
\\[-8pt]
\\[-8pt]
&&\qquad {} + \biggl\llvert \mathcal{E}^{\prime} \biggl[ v \biggl(
\frac
{i-1}{n} , B_n x + B_n Y_{1} \biggr)
\biggr] - v \biggl( \frac
{i-2}{n} , B_n x \biggr) \biggr\rrvert
\nonumber
\\
&&\quad < \frac{2 K_2}{n^{ 1+ \gamma/ \alpha}}+ \frac{\varepsilon}{n}
\nonumber
\end{eqnarray}
for $2 \leq i \leq n-1$, $x \in\mathbb{R}$, and $n \geq N$.

Since $Y_{i+1}$ is independent from $ (Y_1 , \dots, Y_i )$,
repeated application of (\ref{bd4a}) shows that for $n \geq N$,
\begin{equation}
\label{bd4b} \mathcal{E}^{\prime} \Biggl[ \sum
_{i = 2}^{n-1} D^{n}_{i} \Biggr]
< ( n -2 ) \biggl( \frac{2 K_2}{n^{1 +
\gamma/ \alpha}}+ \frac{\varepsilon}{n} \biggr) <
\frac{2
K_2}{n^{\gamma/ \alpha}}+ \varepsilon
\end{equation}
and
\begin{equation}
\label{bd4c} \mathcal{E}^{\prime} \Biggl[ \sum
_{i = 2}^{n-1} D^{n}_{i} \Biggr]
> - ( n -2 ) \biggl( \frac{2 K_2}{n^{1
+ \gamma/ \alpha}}+ \frac{\varepsilon}{n} \biggr) > -
\frac{2
K_2}{n^{\gamma/ \alpha}}- \varepsilon.
\end{equation}

We only need to combine our bounds above and invoke hypothesis (ii) to
finish the proof. By (\ref{bd1}), (\ref{bd2}), (\ref{bd3}), (\ref
{bd4b}) and (\ref{bd4c}),
\begin{eqnarray*}
&&\mathcal{E}^{\prime} \bigl[ v ( 1 , B_n S_n )
\bigr] - v (0, 0 )
\\
&&\quad = \mathcal{E}^{\prime} \Biggl[ v ( 1 , B_n
S_n ) - v \biggl( \frac{n-1}{n} , B_n
S_{n} \biggr) + \sum_{i =
2}^{n-1}C^{n}_{i}
+ \sum_{i = 2}^{n-1}D^{n}_{i}
+ v \biggl( \frac{1}{n} , B_n S_{2} \biggr) - v (0 ,
0 ) \Biggr]
\\
&&\quad \leq\mathcal{E}^{\prime} \biggl[ \biggl\llvert v ( 1 ,
B_n S_n ) - v \biggl( \frac{n-1}{n} ,
B_n S_{n} \biggr) \biggr\rrvert \biggr]+ \sum
_{i = 2}^{n-1} \mathcal{E}^{\prime} \bigl[ \bigl
\llvert C^{n}_{i} \bigr\rrvert \bigr] +
\mathcal{E}^{\prime} \Biggl[ \sum_{i = 2}^{n-1}
D^{n}_{i} \Biggr]
\\
&&\qquad {} + \mathcal{E}^{\prime} \biggl[ \biggl\llvert v \biggl(
\frac
{1}{n} , B_n S_{2} \biggr) - v (0 , 0 ) \biggr
\rrvert \biggr]
\\
&&\quad < \biggl( K_1 \sqrt{\frac{1}{n}} \biggr) + \biggl(
K_2 \biggl[ \biggl(\frac{1}{n} \biggr)^{\gamma/ \alpha} +
B_n^{\gamma} \mathcal {E}^{\prime} \bigl[\llvert
Y_{1} \rrvert ^{\gamma} \bigr] \biggr] \biggr) + \biggl(
\frac{2 K_2}{n^{\gamma/ \alpha}}+ \varepsilon \biggr)
\\
&&\qquad {} + \biggl( K_1 \sqrt{\frac{1}{n}} + 2
\operatorname{Lip} (\psi ) B_n \mathcal{E}^{\prime} \bigl[
\llvert Y_1 \rrvert \bigr] \biggr)
\end{eqnarray*}
and
\begin{eqnarray*}
&&\mathcal{E}^{\prime} \bigl[ v ( 1 , B_n S_n )
\bigr] - v (0, 0 )
\\
&&\quad > - \biggl( K_1 \sqrt{\frac{1}{n}} \biggr) - \biggl(
K_2 \biggl[ \biggl(\frac{1}{n} \biggr)^{\gamma/ \alpha} +
B_n^{\gamma} \mathcal {E}^{\prime} \bigl[\llvert
Y_{1} \rrvert ^{\gamma} \bigr] \biggr] \biggr) - \biggl(
\frac{2 K_2}{n^{\gamma/ \alpha}}+ \varepsilon \biggr)
\\
&&\qquad {} - \biggl( K_1 \sqrt{\frac{1}{n}} + 2
\operatorname{Lip} (\psi ) B_n \mathcal{E}^{\prime} \bigl[
\llvert Y_1 \rrvert \bigr] \biggr)
\end{eqnarray*}
for $n \geq N$. Since $\varepsilon> 0$ is arbitrary and $
\lim_{n \rightarrow\infty} B_n = 0$, we have
\[
\lim_{n \rightarrow\infty} \mathcal{E}^{\prime} \bigl[ v ( 1 ,
B_n S_n ) \bigr] = v (0, 0 ).
\]
\upqed \end{pf}

\section{Examples}\label{sectexmp}

\begin{example}\label{x1examp}
$X_1$ is in its own domain of normal attraction. While this follows
directly from the $\alpha$-stability of $X_1$, we will demonstrate
this using Theorem~\ref{newmr} as well in order to unpack our main result.

Let $\psi\in C_{\mathrm{b.Lip}}  ( \mathbb{R}  )$ and $u$ be
defined by
\[
u (t,x ) = \mathcal{E} \bigl[ \psi ( x + X_t ) \bigr]
\]
on $ [0 , \infty ) \times\mathbb{R}$. If $\tilde{X}_1$
is an independent copy of $X_1$, then
\begin{eqnarray*}
\mathcal{E} \bigl[ \psi ( a X_1 + b \tilde{X}_1 ) \bigr]
&=& \mathcal{E} \bigl[ \mathcal{E} \bigl[ \psi \bigl( a x + \bigl(
b^{\alpha} \bigr)^{\sfrac{1}{\alpha}} \tilde{X}_1 \bigr)
\bigr]_{ x=X_1} \bigr]
\\
&=& \mathcal{E} \bigl[ u \bigl( b^{\alpha} , a X_1 \bigr)
\bigr]
\\
&= &u \bigl( a^{\alpha} + b^{\alpha} , 0 \bigr)
\\
&=& \mathcal{E} \bigl[ \psi \bigl( \bigl(a^{\alpha} + b^{\alpha
}
\bigr)^{\sfrac{1}{\alpha}} X_1 \bigr) \bigr]
\end{eqnarray*}
for any $a$, $b \geq0$ by Lemmas \ref{stabprocprop} and \ref{nnlem51}, that is, $X_1$ is $\alpha$-stable. Exercise 3.20 in \cite
{peng+sc} implies that the same relation actually holds for a broader
class of maps. In particular,
\begin{eqnarray*}
2^{\sfrac{1}{\alpha}} \mathcal{E} [ X_1 ] &=& \mathcal {E} \bigl[
\mathcal{E} [ x + \tilde{X}_1 ]_{ x = X_1} \bigr]
\\
&=& \mathcal{E} \bigl[ X_1 + \mathcal{E} [ X_1 ] \bigr]
\\
&=& 2 \mathcal{E} [ X_1 ] ,
\end{eqnarray*}
so
\[
\mathcal{E} [ X_1 ] = 0.
\]
It follows similarly that
\[
\mathcal{E} [- X_1 ] = 0.
\]

We know
\[
\mathcal{E} \bigl[\llvert X_{1} \rrvert \bigr] < \infty
\]
from Lemma~\ref{absvalprop}.

To check the final hypothesis, let $0 < h < 1$ and $v$ be the unique
viscosity solution of (\ref{vPIDE}). Then for all $ (t , x
 ) \in [ 0 , 1 ] \times\mathbb{R}$,
\begin{eqnarray*}
&&n \biggl\llvert \mathcal{E} \bigl[\delta_{B_n X_1} v ( t , x ) \bigr] -
\biggl( \frac{1}{n} \biggr) \sup_{k_{\pm}
\in K_{\pm}} \biggl\{ \int
_{\mathbb{R}} \delta_z v ( t, x ) F_{k_{\pm}} (dz )
\biggr\} \biggr\rrvert
\\
&&\quad = n \biggl\llvert \mathcal{E} \bigl[ v ( t , x + B_n
X_{1} ) \bigr]- v ( t , x ) + \biggl( \frac{1}{n} \biggr)
\partial_t v ( t , x ) \biggr\rrvert
\\
&&\quad = n \biggl\llvert v \biggl( t - \frac{1}{n} , x \biggr) - v ( t ,
x ) + \biggl( \frac{1}{n} \biggr) \partial_t v ( t , x )
\biggr\rrvert
\\
&&\quad = \biggl\llvert \frac{v  ( t - \sfrac{1}{n} , x  ) - v  ( t
, x  )}{1/n} + \partial_t v ( t , x )
\biggr\rrvert
\\
&&\quad \leq\frac{K_2}{n^{\gamma/ \alpha}}
\end{eqnarray*}
by (\ref{relatvu}), (\ref{reg}) and Lemma~\ref{nnlem51}. Here,
$b=1$ or, equivalently,
\[
B_n = \frac{1}{ n^{1/ \alpha}}.
\]
Abusing notation, Theorem~\ref{newmr} shows that
\[
B_n S_n \xrightarrow{d} X_1
\]
as $n \rightarrow\infty$.
\end{example}

\begin{example}\label{bignlexmp}
Up to some ``uniformity'' assumptions, this example has a
straightforward interpretation.
\begin{quote}
Let the uncertainty subset of distributions (see \cite{peng+sc}) of
$Y_1$ be given by $ \{ \mathbb{P}_\theta: \theta\in\Theta
 \}$. If for all $\theta\in\Theta$, a classical random
variable with distribution $\mathbb{P}_\theta$ is in the domain of
normal attraction of a classical $\alpha$-stable random variable with
triplet $\theta$, then $Y_1$ is in the domain of normal attraction of $X_1$.
\end{quote}

Let $b$, $M >0$ and $f$ be a non-negative function on $\mathbb{N}$
tending to zero as $n \rightarrow\infty$.
For each $k_{\pm} \in K_{\pm}$, let $W_{k_{\pm}}$ be a classical
random variable such that
\begin{enumerate}[(iii)]
\item[(i)]$W_{k_{\pm}}$ has mean zero;
\item[(ii)]$W_{k_{\pm}}$ has a cdf $F_{W_{k_{\pm}}}$ of the form
\begin{equation}
\label{domnorm521} F_{W_{k_{\pm}}} ( z ) = \cases{ \bigl[ b^{\alpha} ( k_- /
\alpha ) + \beta_{1,k_\pm} ( z ) \bigr] \dfrac{1}{\llvert  z \rrvert  ^\alpha}, & $z < 0$,
\cr
1 -
\bigl[ b^{\alpha} ( k_+ / \alpha ) + \beta_{2,k_\pm
} ( z ) \bigr]
\dfrac{1}{z^\alpha}, & $z > 0$, } %
\end{equation}
for some continuously differentiable functions $\beta_{1,k_{\pm}}$ on
$ ( -\infty, 0  ]$ and $\beta_{2,k_{\pm}}$ on $ [ 0
, \infty )$ with
\[
\lim_{z \rightarrow- \infty} \beta_{1,k_{\pm}} ( z ) = \lim
_{z \rightarrow\infty} \beta_{2,k_{\pm}} ( z ) = 0 ;
\\
\]
\item[(iii)] the following quantities are all less than $M$:
\begin{eqnarray*}
&& \biggl\llvert \int_{-\infty}^{-1}
\frac{\beta_{1,k_{\pm}}  ( z
 )}{ ( -z )^{\alpha}}\,dz \biggr\rrvert , \qquad \biggl\llvert \int
_{-1}^{0} \frac{ \beta^{\prime}_{1,k_{\pm}}  ( z
 )}{ ( -z )^{\alpha-1}}\,dz \biggr\rrvert ,
\qquad \int_{-1}^{0} \frac{\llvert   {-} \beta_{1,k_{\pm}}^{\prime}  ( z
 ) z + \alpha\beta_{1,k_{\pm}}  ( z  )\rrvert   }{
 (- z  )^{\alpha- 1}}\,dz ,
\\
&& \biggl\llvert \int_{1}^{\infty} \frac{\beta_{2,k_{\pm}}  ( z
 )}{z^{\alpha}}
\,dz \biggr\rrvert ,\qquad \biggl\llvert \int_{0}^{1}
\frac{ \beta^{\prime}_{2,k_{\pm}}  ( z
)}{z^{\alpha-1}}\,dz \biggr\rrvert ,\qquad \int_{0}^{1}
\frac
{\llvert   {-} \beta_{2,k_{\pm}}^{\prime}  ( z  ) z + \alpha
\beta_{2,k_{\pm}}  ( z  )\rrvert   }{ z^{\alpha- 1}}\,dz;
\end{eqnarray*}
and
\item[(iv)] the following quantities are less than $f  ( n  )$ for
all $n$:
\begin{eqnarray*}
&& \bigl\llvert \beta_{2,k_{\pm}} \bigl( B_{n}^{-1}
\bigr) \bigr\rrvert ,\qquad \int_{1}^{\infty}
\frac{\llvert  \beta_{2,k_{\pm}}  (
B_{n}^{-1} z  ) \rrvert  }{z^{\alpha}}\,dz ,\qquad \int_{0}^{1}
\frac{\llvert   \beta_{2,k_{\pm}}  ( B_{n}^{-1} z
)\rrvert   }{ z^{\alpha- 1}}\,dz
\\
&& \bigl\llvert \beta_{1,k_{\pm}} \bigl( - B_{n}^{-1}
\bigr) \bigr\rrvert ,\qquad \int_{-\infty}^{-1}
\frac{\llvert  \beta_{1,k_{\pm}}
 ( B_{n}^{-1} z  ) \rrvert  }{ (-z )^{\alpha}}\,dz ,\qquad \int_{-1}^{0}
\frac{\llvert   \beta_{1,k_{\pm}}
 ( B_{n}^{-1} z  )\rrvert   }{  (-z )^{\alpha-
1}}\,dz .
\end{eqnarray*}
\end{enumerate}
Note that by (ii) alone, the terms in (iii) are finite and the terms in
(iv) approach zero as \mbox{$n \rightarrow\infty$}. In other words, the
content of (iii) and (iv) is that uniform bounds and minimum rates of
convergence exist.

Define an operator $\mathcal{E}^{\prime}$ on a space $\mathcal{H}$
of suitable functions by
\[
\mathcal{E}^{\prime} [ \varphi ] = \sup_{k_{\pm} \in K_{\pm}} \int
_{\mathbb{R}} \varphi ( z )\,dF_{W_{k_{\pm}}} ( z )
\]
for all $\varphi\in\mathcal{H}$. The exact composition of $\mathcal
{H}$ is irrelevant for our purposes here. Clearly, $ ( \mathbb{R}
, \mathcal{H} , \mathcal{E}^{\prime}  )$ is a sublinear
expectation space.

Let $Y_1$ be the random variable on this space defined by
\[
Y_1 ( x ) = x
\]
for all $x \in\mathbb{R}$. We will use Theorem~\ref{newmr} to show that
\[
B_n S_n \xrightarrow{d} X_1
\]
as $n \rightarrow\infty$. Most of the difficulties have already been
addressed during our discussion of the classical case in Section~\ref{mrsection}.

Since each $W_{k_{\pm}}$ has mean zero,
\[
\mathcal{E}^{\prime} [ Y_1 ] = \sup_{k_{\pm} \in K_{\pm}}
\int_{\mathbb{R}} z \, d F_{W_{k_{\pm}}} ( z ) = 0
\]
and
\[
\mathcal{E}^{\prime} [- Y_1 ] = \sup_{k_{\pm} \in K_{\pm}}
\int_{\mathbb{R}} - z \, d F_{W_{k_{\pm}}} ( z ) = 0.
\]
After recalling that $K_{\pm} \subset ( \lambda, \Lambda
)$, (iii) gives
\begin{eqnarray*}
\mathcal{E}^{\prime} \bigl[ \llvert Y_1 \rrvert \bigr] <
\infty
\end{eqnarray*}
using (\ref{absvalint1}) and (\ref{domnorm521}). Observe that
we are solving (\ref{domnorm521}) for the obvious expressions to
obtain uniform bounds on the terms
\[
\bigl\llvert \beta_{2,k_{\pm}} ( 1 ) \bigr\rrvert ,\qquad \bigl\llvert
\beta_{1,k_{\pm}} ( - 1 ) \bigr\rrvert ,\qquad \biggl\llvert \int
_{0}^{1} \frac{b^\alpha k_{+} + \alpha\beta
_{2,k_{\pm}}  ( z  )}{z^{\alpha}}\,dz \biggr\rrvert
\]
and
\[
\biggl\llvert \int_{-1}^{0} \frac{b^\alpha k_{-} + \alpha\beta_{1,k_{\pm
}}  ( z  )}{ ( -z )^{\alpha}}
\,dz \biggr\rrvert .
\]

To check the remaining hypothesis, let $0 < h < 1$, $\psi\in C_{\mathrm{b.Lip}}
 ( \mathbb{R}  )$, and $v$ be the unique viscosity solution
of (\ref{vPIDE}). The techniques of (\ref{linequivdisthyp2})
demonstrate that
\begin{eqnarray*}
&&n \biggl\llvert \mathcal{E}^{\prime} \bigl[ \delta_{B_n Y_1} v (
t, x ) \bigr] - \biggl( \frac{1}{n} \biggr) \sup_{k_{\pm} \in K_{\pm}}
\biggl\{ \int_{\mathbb{R}} \delta_z v ( t, x )
F_{k_{\pm}} (dz ) \biggr\} \biggr\rrvert
\\
&&\quad \leq \biggl( \frac{1}{b^{\alpha}} \biggr) \sup_{k_{\pm} \in K_{\pm}}
\biggl\llvert \int_{\mathbb{R}} \delta_z v ( t, x )
\biggl( \frac{ \beta_{1,k_{\pm}}^{\prime}  (
B_{n}^{-1} z  ) \llvert  B_{n}^{-1} z \rrvert   + \alpha\beta
_{1,k_{\pm}}  ( B_{n}^{-1} z  ) }{ \llvert   z \rrvert  ^{\alpha+ 1}} \mathbf{1}_{ ( -\infty, 0  )} ( z )
\\
&&\qquad {} + \frac{- \beta
_{2,k_{\pm}}^{\prime}  ( B_{n}^{-1} z  ) \llvert  B_{n}^{-1}
z \rrvert   + \alpha\beta_{2,k_{\pm}}  ( B_{n}^{-1} z  )
}{ \llvert   z \rrvert  ^{\alpha+ 1}} \mathbf{1}_{ ( 0 , \infty
)} ( z ) \biggr)\,dz \biggr
\rrvert
\end{eqnarray*}
for $ ( t , x  ) \in [0, 1 ] \times\mathbb{R}$
and $n \geq1$. Combining (\ref{intbd1}), (\ref{intbd2}), (\ref
{intbd3}) and (\ref{intbd4}) with (iii) and (iv) proves that this
last expression approaches zero in the required way.
\end{example}

\begin{appendix}

\section{Interior regularity theory
background}\label{PIDEbackgroundapp}

Interior regularity theory for fully nonlinear integro-differential
equations is rich and well developed. Before describing the
results that we need for our proof, we provide a short discussion of
the literature. Readers new to this field are encouraged to
first consult \cite{wiki} for an introduction.

Some results and methods from the interior regularity theory
for PDEs can be imported to the non-local case after minor
modifications. For other aspects of the theory, this is false. As
described in Section~2 of \cite{sch+rang+kass}, a H\"{o}lder estimate
and the Harnack inequality appear together in the local setting;
however, there are non-local equations for which a H\"{o}lder estimate
holds in the absence of the Harnack inequality. A partial list of other
ways that non-local results can significantly differ from their local
counterparts can be found in \cite{wiki}.

Early work on the regularity of integro-differential equations
focused on equations in divergence form. A survey of these results is
contained in \cite{kass+sch}. For equations in non-divergence form,
\cite{bass+levin} contains the first Harnack inequality and H\"{o}lder
estimate. The equations studied in \cite{bass+levin} are of the form
\[
 \int_{\mathbb{R}^d} \bigl[ w ( x + z ) - w ( x ) - z \nabla w
( x ) {\bf 1}_{B_1} ( z ) \bigr] k ( x, z )\,dz = 0,
\]
where $k$ is a kernel such that
\begin{equation}
\label{symmkern}  k ( x , z ) = k ( x , -z )
\end{equation}
and
\begin{equation}
\label{kernbds}  \frac{\lambda_1}{ \llvert   z \rrvert  ^{d +
\alpha_1}} \leq k ( x, z ) \leq
\frac
{\Lambda_1}{ \llvert   z \rrvert  ^{d + \alpha_1} }
\end{equation}
for some constants $\lambda_1$, $\Lambda_1 > 0$ and $\alpha
_1 \in ( 0 ,2  )$. For a review of the extensions of this
initial work, see \cite{kass+sch}.

The H\"{o}lder estimate in \cite{bass+levin} blows up as
$\alpha_1 \rightarrow2$. Many other early estimates share this
feature. The first paper to prove a H\"{o}lder estimate and Harnack
inequality without this property is \cite{caff+silvestre}. The
equations are of the form
\begin{equation}
\label{caffsilveqn}  \inf_r \sup_s
\biggl\{ \int_{\mathbb{R}^d} \bigl[ w ( x + z ) - w ( x ) -z \nabla w
( x ) \mathbf{1}_{B_1} ( z ) \bigr] k^{rs} ( z )\,dz \biggr\} =
0
\end{equation}
for kernels $k^{rs}$ depending only on $z$ and satisfying (\ref
{symmkern}), (\ref{kernbds}) and an additional smoothness condition.
More precisely, for some fixed positive constants $\rho$ and $C$,
\[
 \int_{\mathbb{R}^d \setminus B_{\rho}} \frac{ \llvert   k  ( z  ) - k  ( z-
\varepsilon ) \rrvert  }{\llvert  \varepsilon\rrvert  }\,dz \leq C
\]
whenever
\[
 \llvert \varepsilon\rrvert < \frac{\rho}{2}.
\]
The paper culminates in a $C^{1,\gamma}$ estimate for the
solution of (\ref{caffsilveqn}).

These findings have been extended in a number of ways. For
instance, references such as \cite{silv,schwab+silv,MR2995098,lara+dav+crit+dr} and \cite{lara+dav+sig+al} study
equations with non-symmetric kernels, that is, kernels that do not
satisfy (\ref{symmkern}). Other examples of recent work include \cite
{caff+silv+evans,serra+classical} and \cite{kriventsov}.

We now collect the definitions and
results from \cite{lara+dav+crit+dr} and \cite{lara+dav+sig+al}
that we need for our proof. These references describe
properties of the solutions to a broad class of non-local fully
nonlinear parabolic equations of the form
\[
\partial_t w ( t , x ) - I w ( t, x ) = f ( t ).
\]
%
Due to the general nature of these equations, \cite
{lara+dav+crit+dr} and \cite{lara+dav+sig+al} are quite technical.
Since (\ref{uPIDEback}) is an easy case of the equations studied in
these papers, we will simplify
this material and present only the version that we need for our argument.

\begin{notn}
Let
\[
\mathfrak{C} _{\tau, r} ( t, x ) := ( t - \tau, t ] \times
(x -r ,x+ r ).
\]
We write $\mathfrak{C}_{\tau,r}$ for the
cylinder $\mathfrak{C}_{ \tau,r}  ( 0 , 0  )$.
For suitable functions $w$, let
\begin{eqnarray*}
 \tilde{\delta}_z w ( t , x ) &:=& w ( t , x + z ) - w ( t , x ) -
\partial_x w ( t, x ) \mathbf{1}_{ ( -1, 1 )} ( z ) z ;
\\
\llVert w \rrVert _{L^1  ( \nu )} &:=& \int_{\mathbb{R}} \bigl
\llvert w ( z ) \bigr\rrvert \min \bigl( 1 , \llvert z \rrvert ^{-1-\alpha}
\bigr)\,dz ; \quad \mbox{and}
\\
{} [ w ]_{C^{0,1}  (  (t_0, t_1  ] \mapsto
L^1  ( \nu )  ) } &:=& \sup_{  ( t
- \tau, t  ] \subseteq ( t_0 , t_1  ] } \frac{
\llVert   w  ( t , \cdot ) - w  ( t - \tau, \cdot
 ) \rrVert  _{L^1  ( \nu )}}{\tau}.
\end{eqnarray*}
We also let
\[
b_{k_{\pm}} :=  ( k_- - k_+ )
\int
_1^{\infty} \frac{dz }{ z^{\alpha} }
\]
for all $k_{\pm} \in K_{\pm}$.
\end{notn}

In the literature, one also works frequently with cylinders of
the form
\[
 ( t - \tau^\alpha, t ] \times (x -r ,x+ r )
\]
due to their convenient scaling properties. We introduce
\[
\llVert \cdot \rrVert _{L^1  ( \nu )}
\]
and
\[
 [ \cdot ]_{C^{0,1}  (  (t_0, t_1
 ] \mapsto L^1  ( \nu )  ) }
\]
due to their role in upcoming H\"{o}lder estimates, namely,
Lemmas \ref{inithold} and \ref{powerest}. The symbols $\tilde
{\delta}_z$ and $b_{k_{\pm}}$ facilitate the identification of (\ref
{uPIDEback}) with the equations studied \cite{lara+dav+crit+dr} and
\cite{lara+dav+sig+al}. Observe that for all $k_{\pm} \in K_{\pm}$
and suitable functions $w$,
\begin{eqnarray}
\label{deltildvsdel}  \int_{\mathbb{R}} \delta_z w (
t , x ) F_{k_{\pm}} ( dz ) = b_{k_{\pm}} \partial_x w ( t
, x ) + \int_{\mathbb{R}} \tilde {\delta}_z w ( t , x )
F_{k_{\pm}} ( dz ) .
\end{eqnarray}

\begin{defn}
Since $K_{\pm} \subset ( \lambda, \Lambda )$, we can
pick $\beta> 0$ such that
\[
\sup_{k_{\pm} \in K_{\pm}} \biggl\{ \sup_{r \in (0 , 1  )} \biggl\{
r^{\alpha- 1} \biggl\llvert b_{k_\pm
} + \int_{  ( - 1, 1 ) \backslash ( -r
, r  )}
z F_{k_{\pm}} ( dz ) \biggr\rrvert \biggr\} \biggr\} \leq\beta.
\]
Let $\mathcal{L}_0$ be the family of operators
\[
w (t, x ) \mapsto b \partial_x w ( t , x ) + \int
_{\mathbb{R}} \tilde{\delta}_z w ( t , x )
\frac{ k  (z  )}{\llvert   z \rrvert  ^{1 + \alpha
}}\,dz ,
\]
where $k$ is a kernel and $b$ is a constant such that $\lambda\leq k
\leq\Lambda$ and
\[
\sup_{r \in (0 , 1  )} r^{\alpha- 1} \biggl\llvert b + \int
_{ ( - 1, 1 ) \backslash ( -r
, r  )} \frac{z k  (z )}{\llvert  z \rrvert  ^{1+\alpha
}}\,dz \biggr\rrvert \leq\beta.
\]
We say that an operator in $\mathcal{L}_0$ is in $\mathcal{L}_1$ if
\[
\bigl\llvert \partial_z k ( z ) \bigr\rrvert \leq
\frac{\Lambda
}{\llvert   z \rrvert  },
\]
and an operator in $\mathcal{L}_1$ is in $\mathcal{L}_2$ if
\[
\bigl\llvert \partial_{zz}^2 k ( z ) \bigr\rrvert \leq
\frac
{\Lambda}{ \llvert   z \rrvert  ^2}.
\]
\end{defn}

The stronger regularity requirements on the kernels (in
$\mathcal{L}_2$, say, compared to those in $\mathcal{L}_0$) give rise
to stronger regularity results.  All of the operators
\[
w (t, x ) \mapsto b_{k_{\pm}} \partial_x w ( t , x ) + \int
_{\mathbb{R}} \tilde{\delta}_z w ( t , x )
F_{k_{\pm}} ( dz )
\]
are in each of these families. As we will soon see in (\ref{ubareqn}), we will be especially interested in the operator $I$
defined by
\[
I w ( t, x ) = \inf_{k_{\pm} \in K_{\pm}} \biggl\{ b_{k_{\pm}}
\partial_x w ( t , x ) + \int_{\mathbb{R}} \tilde{
\delta}_z w ( t , x ) F_{k_{\pm}} ( dz ) \biggr\}.
\]
$I$ is a specific case of an {\it extremal operator}.

\begin{defn}
For a collection of operators $\mathcal{L} \subseteq\mathcal{L}_0$,
define the {\it extremal operators} $\mathcal{M}^{+}_{\mathcal{L}}$
and $\mathcal{M}^{-}_{\mathcal{L}}$ by
\[
\mathcal{M}^{+}_{\mathcal{L}} = \sup_{L \in\mathcal
{L}} L
\quad \mbox{and}\quad  \mathcal{M}^{-}_{\mathcal{L}} = \inf_{L \in\mathcal{L}}
L.
\]
\end{defn}
$I$ has a number of other key properties including the
following.\footnote{Though we will not emphasize this point,
we remark in passing that $I w ( t, x  ) $ is well-defined
for any $w ( t, \cdot ) \in C^{1,1}  ( x  ) \cap
L^1  ( \nu )$ (see Section~2 of \cite{lara+dav+sig+al}).}
\begin{enumerate}[(iii)]
\item[(i)]$I 0 = 0$.
\item[(ii)]$I$ is {\it uniformly elliptic} with respect to $\mathcal{L}_j$,
that is,
\[
\mathcal{M}^{-}_{\mathcal{L}_j} ( w_1 - w_2 )
\leq I w_1 - I w_2 \leq\mathcal{M}^{+}_{\mathcal{L}_j}
( w_1 - w_2 ).
\]
\item[(iii)]$I$ is {\it translation invariant}, that is,
\[
I \bigl( w ( t_0+ \cdot, x_0 + \cdot ) \bigr) ( t, x ) =
( I w ) (t_0 + t, x_0 + x ).
\]
\end{enumerate}




(i) is trivial. See Section~2 of \cite{lara+dav+sig+al} for
(ii). Since $I$ has constant coefficients, we get (iii). We highlight
these classes of operators and properties of $I$ for the convenience of
the reader comparing the next three results to their original versions
(see Theorem~2.3 in \cite{lara+dav+sig+al} for Lemma~\ref{inithold};
Theorems 1.1, 2.4, and 2.5 in \cite{lara+dav+sig+al} for Lemma~\ref
{powerest}; and Theorem~3.3 in \cite{lara+dav+crit+dr} for Lemma~\ref
{diffviscsol}).\footnote{A number of related results exist
in the literature. We mention only a small sample. Theorem~12.1 in
\cite{caff+silvestre}, Theorem~1.1 in \cite{sch+rang+kass} and
Theorem~7.1 in \cite{schwab+silv} are $C^\gamma$ estimates along the
lines of Lemma~\ref{inithold}. Theorem~8.1 in \cite{schwab+silv},
Theorem~13.1 in \cite{caff+silvestre}, Theorem~1.1 in \cite
{caff+silv+evans} and Theorem~1.1 in \cite{serra+classical} contain
$C^{1,\gamma}$ or $C^{\alpha+\gamma}$ estimates similar to those in
Lemma~\ref{powerest}. Like Lemma~\ref{diffviscsol}, Theorem~5.9 in
\cite{caff+silvestre} and Lemma~3.2 in \cite{Silvestre20112020}
investigate the difference of viscosity solutions.}

\begin{lem}\label{inithold}
Let $w$ satisfy
\begin{eqnarray*}
\partial_t w - M^{+}_{\mathcal{L}_0} w&\leq&0,
\\
\partial_t w - M^{-}_{\mathcal{L}_0} w &\geq&0
\end{eqnarray*}
in the viscosity sense on $\mathfrak{C}_{1,1}$. There
is some $\gamma\in ( 0, 1  )$ and $C > 0$ depending only
on $\lambda$, $\Lambda$, and $\beta$ such that for every $ (
t_0, x_0  )$, $ (t_1, x_1  ) \in\mathfrak{C}_{1/2 , 1/2}$,
\[
\frac{ \llvert   w  ( t_0, x_0  ) - w  (t_1, x_1
) \rrvert   }{ ( \llvert   t_0 - t_1 \rrvert  ^{1/\alpha} + \llvert
x_0 - x_1 \rrvert    )^{\gamma}} \leq C \llVert w \rrVert _{L^1
 (  (-1, 0  ] \mapsto L^1  ( \nu )  )}.
\]
\end{lem}

\begin{lem}\label{powerest}
Let $w$ satisfy
\[
\partial_t w - Iw = 0
\]
in the viscosity sense on $\mathfrak{C}_{1,1}$ . There
is some $\gamma\in ( 0, 1  )$ and $C > 0$ depending only
on $\lambda$, $\Lambda$, and $\beta$ such that for every $ (
t_0, x_0  )$, $ (t_1, x_1  ) \in\mathfrak{C}_{1/2 , 1/2}$,
\[
\bigl\llvert \partial_x w ( t_0, x_0 )
\bigr\rrvert + \frac{
\llvert   \partial_x w ( t_0, x_0  ) - \partial_x w  (t_1, x_1  )\rrvert  }{ ( \llvert   t_0 - t_1 \rrvert  ^{1/\alpha}+ \llvert   x_0 - x_1 \rrvert    )^{\gamma}} \leq C \llVert w \rrVert _{L^1  (  (-1, 0  ] \mapsto L^1
 ( \nu )  )}
\]
and
\[
\bigl\llvert \partial_t w ( t_0, x_0 )
\bigr\rrvert + \frac{
\llvert   \partial_t w  ( t_0, x_0  ) - \partial_t w  (t_1, x_1  )\rrvert  }{ ( \llvert   t_0 - t_1 \rrvert  ^{1/\alpha} + \llvert   x_0 - x_1 \rrvert    )^{\gamma}} \leq C [ w ]_{C^{0,1}  (  (-1, 0  ] \mapsto L^1
 ( \nu )  ) } .
\]
We also have
\[
\llVert w \rrVert _{C^{\alpha+ \gamma}  ( \mathfrak{C}_{1/2,1/2}  )} \leq C \bigl( \llVert w \rrVert
_{L^1
 (  (-1, 0  ] \mapsto L^1  ( \nu )
)} + [ w \mathbf{1}_{ (-1,1 )^c} ]_{C^{0,1}  (
 (-1, 0  ] \mapsto L^1  ( \nu )  ) } \bigr).
\]
\end{lem}

\begin{lem}\label{diffviscsol}
Let $w_1$, $w_2$ satisfy
\[
\partial_t w_i - I w_i = 0
\]
in the viscosity sense on some domain $\Omega$. Then
\begin{eqnarray*}
\partial_t ( w_1 - w_2 ) -
M^{+}_{\mathcal{L}_0} ( w_1 - w_2 ) &\leq&0,
\\
\partial_t ( w_1 - w_2 ) -
M^{-}_{\mathcal{L}_0} ( w_1 - w_2 ) &\geq&0
\end{eqnarray*}
also holds in the viscosity sense on $\Omega$.
\end{lem}

We will need one more result (for the original version, see Lemma~5.6
and the proof of Corollary~5.7 in \cite{caff+cabre}). It is the key to
a standard technique from the literature allowing one to repeatedly
apply an estimate such as Lemma~\ref{inithold} in order to obtain a
higher regularity estimate.

\begin{lem}\label{lem56cc}
Let $0< \beta_1 \leq1$, $0 < \beta_2 < 1$, $L >0$, and $w \in
L^{\infty}  (  [-1, 1  ]  )$ satisfy
\[
\llVert w \rrVert _{ L^{\infty}  (  [-1, 1  ]
)} \leq L.
\]
For $0 < \llvert   h_0 \rrvert   \leq1$, define $w_{\beta_1, h_0}$ by
\[
w_{\beta_1, h_0} ( x ) = \frac{w  (x
+ h_0  ) - w  ( x )}{\llvert  h_0\rrvert   ^{\beta_1}}
\]
for all $x \in I_{h_0}$, where $I_{h_0} =  [ -1 , 1 - h_0
]$ if $h_0 >0$ and $I_{h_0} =  [ -1 -h_0 , 1  ]$ if $h_0 <
0$. Suppose that
\[
w_{\beta_1, h_0} \in C^{\beta_2} ( I_{h_0} )
\]
and
\[
\llVert w_{\beta_1, h_0} \rrVert _{ C^{\beta_2}  ( I_{h_0}
 ) } \leq L
\]
for any $0 < \llvert   h_0 \rrvert   \leq1$.
\begin{enumerate}[(iii)]
\item[(i)] If $\beta_1 + \beta_2 < 1$, then
\[
w \in C^{\beta_1 + \beta_2} \bigl( [-1 , 1 ] \bigr)
\]
and
\[
\llVert w \rrVert _{ C^{\beta_1 + \beta_2}  (  [-1 ,1
 ]  ) } \leq C L.
\]
\item[(ii)] If $\beta_1 + \beta_2 > 1$ and $\beta_1 \ne1$, then
\[
w \in C^{0,1} \bigl( [-1 , 1 ] \bigr)
\]
and
\[
\llVert w \rrVert _{ C^{0,1}  (  [-1 ,1  ]  ) } \leq C L.
\]
\item[(iii)] If $\beta_1 = 1$, then $w \in C^{1,\beta_2}  (  [-1 ,
1  ]  )$ and
\[
\llVert w \rrVert _{ C^{1,\beta_2}  (  [-1 ,1  ]
 ) } \leq C L.
\]
\end{enumerate}
In any of these cases, $C$ depends only on $\beta_1 + \beta_2$.
\end{lem}

We will often apply these results on different domains than we have
listed above without comment. For instance, we might use Lemma~\ref
{powerest} on $\mathfrak{C}_{1,1}  ( t, x
)$ or Lemma~\ref{lem56cc} on an arbitrary closed interval. These
``new'' results are obtained merely by translating or rescaling, both
standard routines in the literature. As an example of such an
operation, notice that if $w$ satisfies
\[
\partial_t w - Iw = 0
\]
in the viscosity sense on $ ( t_1 , t_2  ] \times\Omega$,
then $\tilde{w}$ defined by
\[
\tilde{w} ( t, x ) = w \bigl( r^{\alpha} t + t_0 , rx +
x_0 \bigr)
\]
satisfies
\[
\partial_t \tilde{w} - I \tilde{w} = 0
\]
in the viscosity sense on
\[
\bigg( \frac{t_1-t_0}{r^{\alpha}} , \frac{t_2-t_0}{r^{\alpha}} \bigg] \times\frac{\Omega- x_0}{r}
\]
(see Section~2.1.1 of \cite{lara+dav+crit+dr}). Further information
can be found in \cite{lara+dav+crit+dr,lara+dav+sig+al} and
\cite{caff+cabre}.

\section{Proof of Proposition \texorpdfstring{\protect\ref{uPIDEregprop}}{2.10}}\label{PIDEappend}

In the hope of keeping the number of constants in our argument at a
reasonable level, we will not issue a new subscript each time we
introduce a new constant $B$ below. Also, we will write $\bar{u}$
instead of $-u$. From (\ref{uPIDEback}) and (\ref{deltildvsdel}), $\bar{u}$ is a viscosity solution of
\begin{eqnarray}
\label{ubareqn} \partial_t \bar{u} ( t, x ) - I \bar{u} ( t, x ) &=&
0 ,\qquad  (t, x ) \in (0 , \infty ) \times\mathbb{R}
\nonumber
\\[-8pt]\\[-8pt]
\bar{u} ( 0 ,x ) &=& -\psi ( x ),\qquad  x \in \mathbb{R}.\nonumber
\end{eqnarray}
It suffices to show that parts (i)--(iv) of Proposition~\ref{uPIDEregprop} hold for $\bar{u}$ and (\ref{ubareqn}).

The quantities
\[
[ \bar{u} ]_{C^{0,1}  (  (t_0, t_1  ]
\mapsto L^1  ( \nu )  ) }
\]
play a crucial role in Lemma~\ref{powerest}, so our first goal will
be to control them for $t_0$ greater than some positive number. We will
do this by showing that $\bar{u}$ is uniformly Lipschitz as a function
of time for times above some lower bound. Achieving a Lipschitz
estimate can be done using a standard strategy. Specifically, we will
begin by obtaining an initial $C^{\gamma/ \alpha}$ estimate from
Lemma~\ref{inithold}. Lemma~\ref{diffviscsol} will allow us to
apply Lemma~\ref{inithold} to get a $C^{\gamma/ \alpha}$ estimate
for the incremental quotients of $\bar{u}$. Then Lemma~\ref{lem56cc} will give that $\bar{u}$ is $C^{2\gamma/ \alpha}$ in time. We
will repeat these steps to show that $\bar{u}$ is $C^{3\gamma/ \alpha
}$ in time, $C^{4\gamma/ \alpha}$ in time, and so on until we
conclude that $\bar{u}$ is Lipschitz in time.

Since
\[
\mathcal{M}_{\mathcal{L}_0}^{-} w \leq I w \leq\mathcal
{M}_{\mathcal{L}_0}^{+} w ,
\]
$\bar{u}$ satisfies
\begin{eqnarray*}
\partial_t \bar{u} - M^{+}_{\mathcal{L}_0} \bar{u} &
\leq&0,
\\
\partial_t \bar{u} - M^{-}_{\mathcal{L}_0} \bar{u} &
\geq&0
\end{eqnarray*}
in the viscosity sense on $ (0 , \infty ) \times\mathbb
{R}$. For any $\bar{t} > 1$,
\begin{eqnarray*}
\bigl\llVert \bar{u} ( \bar{t} + \cdot, \cdot ) \bigr\rrVert _{L^1  (  (-1, 0  ] \mapsto L^1  ( \nu )
 )} &=&
\int_{-1}^0 \int_{\mathbb
{R}} \bigl
\llvert \bar{u} ( \bar{t} + t , z ) \bigr\rrvert \min \bigl( 1 , \llvert z \rrvert
^{-1-\alpha} \bigr)\,dz \,dt
\\
&\leq&\llVert \psi\rrVert _{L^{\infty}  ( \mathbb{R}
)} \int_{-1}^0
\int_{\mathbb{R}} \min \bigl( 1 , \llvert z \rrvert ^{-1-\alpha}
\bigr)\,dz \,dt
\end{eqnarray*}
by Lemma~\ref{nnlem53}. Lemma~\ref{inithold} implies that for some
$B$, $\gamma> 0$,
\begin{equation}
\label{initholdpf} \frac{ \llvert   \bar{u}  ( t_0, x_0  ) - \bar{u}
 (t_1, x_1  ) \rrvert   }{ ( \llvert   t_0 - t_1 \rrvert  ^{1/\alpha} + \llvert
  x_0 - x_1 \rrvert    )^{\gamma} } \leq B
\end{equation}
for every $ ( t_0, x_0  )$, $ (t_1, x_1  ) \in
\mathfrak{C}_{1/2 , 1/2}  ( \bar{t} , \bar{x}
 )$ with $\bar{t} > 1$.

For $0 < \llvert   h_0 \rrvert   < 1/2$, define $\bar{u}_{\gamma/ \alpha
, h_0}$ by
\[
\bar{u}_{\gamma/\alpha, h_0} ( t, x ) = \frac{ \bar{u}  (t + h_0, x  ) - \bar{u}  (t ,
x )}{\llvert  h_0\rrvert   ^{\gamma/ \alpha}}
\]
for all $ ( t , x  ) \in [ 1/2, \infty ) \times
\mathbb{R}$. Then
\[
\llVert \bar{u}_{\gamma/\alpha, h_0} \rrVert _{L^{\infty}  (
 ( 1 , \infty ) \times\mathbb{R}  )} \leq B
\]
by (\ref{initholdpf}). Hence,
\begin{eqnarray*}
\bigl\llVert \bar{u}_{\gamma/\alpha, h_0} ( \bar{t} + \cdot, \cdot ) \bigr\rrVert
_{L^1  (  (-1, 0  ] \mapsto L^1
 ( \nu )  )} \leq B \int_{\mathbb{R}} \min \bigl( 1 , \llvert z
\rrvert ^{-1-\alpha} \bigr)\,dz
\end{eqnarray*}
for any $\bar{t} > 2$.

Notice that
\[
\partial_t \bar{u} ( \cdot+ h_0, \cdot ) - I \bar{u} (
\cdot+ h_0, \cdot ) = 0
\]
in the viscosity sense on $ ( 1/2 , \infty ) \times\mathbb
{R}$ because (\ref{ubareqn}) has constant coefficients. Lemma~\ref
{diffviscsol} implies that
\begin{eqnarray*}
\partial_t \bar{u}_{\gamma/ \alpha, h_0} - M^{+}_{\mathcal{L}_0}
\bar{u}_{\gamma/ \alpha, h_0} &\leq&0,
\\
\partial_t \bar{u}_{\gamma/ \alpha, h_0} - M^{-}_{\mathcal{L}_0}
\bar{u}_{\gamma/ \alpha, h_0} &\geq&0
\end{eqnarray*}
in the viscosity sense on $ ( 1/2, \infty ) \times\mathbb
{R}$. For some $B$,
\[
\frac{ \llvert   \bar{u}_{\gamma/ \alpha, h_0} ( t_0, x_0
) - \bar{u}_{\gamma/ \alpha, h_0}  (t_1, x_1  ) \rrvert
}{ ( \llvert   t_0 - t_1 \rrvert  ^{1/\alpha} + \llvert   x_0 - x_1
\rrvert    )^{\gamma} } \leq B
\]
for every $ ( t_0, x_0  )$, $ (t_1, x_1  ) \in
\mathfrak{C}_{1/2 , 1/2}  ( \bar{t} , \bar{x}
 )$ with $\bar{t} > 2$ by Lemma~\ref{inithold}.

Lemma~\ref{lem56cc} shows that for a small $r_1$ (less than $1/4$), we
can find $B$ such that
\[
\bar{u} (\cdot, \bar{x} ) \in C^{2 \gamma/ \alpha} \bigl( [ \bar{t} -
r_1 , \bar{t} + r_1 ] \bigr)
\]
and
\begin{equation}
\label{2gamhold} \bigl\llVert \bar{u} (\cdot, \bar{x} ) \bigr\rrVert
_{C^{2\gamma
/ \alpha}  (  [ \bar{t} - r_1 , \bar{t} + r_1
] )} \leq B
\end{equation}
for $\bar{t} > 2$.

Due to Lemma~\ref{lem56cc}, assume without loss of generality that
$\alpha/ \gamma$ is not an integer. Starting from the incremental quotient
\[
\frac{ \bar{u}  (t + h_0, x  ) - \bar{u}
 (t , x )}{\llvert  h_0\rrvert   ^{2\gamma/ \alpha}},
\]
we can use these steps to produce a $C^{3 \gamma/ \alpha}$ estimate
for $\bar{u}$ in time. By continuing to repeat this procedure, we will
obtain a $C^{4 \gamma/ \alpha}$ estimate, a $C^{5 \gamma/ \alpha}$
estimate, and so on until we obtain a Lipschitz estimate for $\bar{u}$
in time. More precisely, we will find $B$ and a small $r_n$ such that
\[
\bar{u} ( \cdot, \bar{x} ) \in C^{0,1} \bigl( [ \bar{t} -
r_n , \bar{t} + r_n ] \bigr)
\]
and
\[
\bigl\llVert \bar{u} ( \cdot, \bar{x} ) \bigr\rrVert _{C^{0,1}
 (  [ \bar{t} - r_n , \bar{t}+ r_n  ]  )} \leq B
\]
for $\bar{t} > \lceil\alpha/ \gamma\rceil$.

For $t_0$, $t_1 > \lceil\alpha/ \gamma\rceil$,
\begin{eqnarray*}
\bigl\llvert \bar{u} ( t_0 , x_0 ) - \bar{u} (
t_1 , x_0 ) \bigr\rrvert &\leq& \bigl\llvert \bar{u} (
s_0 , x_0 ) - \bar{u} ( s_1 , x_0
) \bigr\rrvert + \cdots+ \bigl\llvert \bar {u} ( s_{N-1} ,
x_0 ) - \bar{u} ( s_{N} , x_0 ) \bigr\rrvert
\\
&\leq& B \llvert s_0 - s_1 \rrvert + \cdots+ B \llvert
s_{N-1} - s_N \rrvert
\\
&=& B \llvert t_0 - t_1 \rrvert,
\end{eqnarray*}
where $t_0 = s_0$, $t_1=s_N$, and $s_{i} < s_{i+1} \leq s_i + 2 r_n$.
This indicates that
\[
\bar{u} ( \cdot, \bar{x} ) \in C^{0,1} \bigl( \bigl( \lceil\alpha/
\gamma\rceil, \infty \bigr) \bigr)
\]
and
\[
\bigl\llVert \bar{u} ( \cdot, \bar{x} ) \bigr\rrVert _{C^{0,1}
 (  ( \lceil\alpha/ \gamma\rceil , \infty )
)} \leq B .
\]

Then $t_0$, $t_1 > \lceil\alpha/ \gamma\rceil$ implies
\begin{eqnarray*}
[ \bar{u} \mathbf{1}_{  ( -1 , 1  )^c} ]_{C^{0,1}
 (  (t_0, t_1  ] \mapsto L^1  ( \nu )
 ) } &\leq& [ \bar{u}
]_{C^{0,1}  (  (t_0,
t_1  ] \mapsto L^1  ( \nu )  ) }
\\
&=& \sup_{  ( t - \tau, t  ] \subseteq ( t_0 , t_1  ] } \frac{ \llVert   \bar{u} ( t , \cdot
 ) - \bar{u}  ( t - \tau, \cdot ) \rrVert  _{L^1
 ( \nu )}}{\tau}
\\
&\leq &B \int_{\mathbb{R}} \min \bigl( 1 , \llvert z \rrvert
^{-1-\alpha} \bigr)\,dz .
\end{eqnarray*}
Lemma~\ref{powerest} gives that for $\bar{t} > \lceil\alpha/
\gamma\rceil$,
\begin{equation}
\label{xderiv} \bigl\llvert \partial_x \bar{u} ( t_0,
x_0 ) \bigr\rrvert + \frac
{ \llvert   \partial_x \bar{u}  ( t_0, x_0  ) - \partial_x
\bar{u}  (t_1, x_1  )\rrvert  }{ ( \llvert   t_0 - t_1
\rrvert  ^{1/\alpha} + \llvert   x_0 - x_1 \rrvert    )^{\gamma} } \leq B
\end{equation}
and
\begin{equation}
\label{tderiv} \bigl\llvert \partial_t \bar{u} ( t_0,
x_0 ) \bigr\rrvert + \frac
{ \llvert   \partial_t \bar{u}  ( t_0, x_0  ) - \partial_t
\bar{u}  (t_1, x_1  )\rrvert  }{ ( \llvert   t_0 - t_1
\rrvert  ^{1/\alpha} + \llvert   x_0 - x_1 \rrvert    )^{\gamma} } \leq B
\end{equation}
for every $ ( t_0, x_0  )$, $ (t_1, x_1  ) \in
\mathfrak{C}_{1/2 , 1/2}  ( \bar{t} , \bar{x}
 )$. It also shows that
\begin{equation}
\label{class} \llVert \bar{u} \rrVert _{C^{\alpha+ \gamma}  ( \mathfrak{C}_{1/2,1/2}  ( \bar{t} , \bar{x}  )
 )} \leq B.
\end{equation}

After suitably rescaling, we see that these inequalities actually hold
for $\bar{t} > (1+h)/2$. Part (i) of Proposition~\ref{uPIDEregprop} then follows from (\ref{xderiv}) and (\ref{tderiv}), while
part (iii) follows from (\ref{class}). From (\ref{tderiv}) and a
simple covering argument, we know that as long as the distance between
$x_0$ and $x_1$ is under some arbitrary bound, we can find $B$ such that
\[
\bigl\llvert \partial_t \bar{u} ( t , x_0 ) -
\partial_t \bar{u} ( t , x_1 ) \bigr\rrvert \leq B
\llvert x_0 - x_1 \rrvert ^{\gamma}
\]
for $t \in [ h , h +1  ]$. Since $\partial_t \bar{u}$ is
bounded on $ [ h , h +1  ] \times\mathbb{R}$, we can drop
the distance constraint and get the second inequality in part (ii). A
similar covering argument finishes the proof of the first inequality
and yields part (ii) of Proposition~\ref{uPIDEregprop}.

It remains to prove part (iv). In this case, the equation for $\bar
{u}$ is
\begin{eqnarray}
\label{uPIDElin} \partial_t \bar{u} (t, x ) - b_{k_{\pm}}
\partial_x \bar{u} ( t , x ) - \int_{\mathbb{R}} \tilde{
\delta}_z \bar{u} ( t , x ) F_{k_{\pm}} ( dz ) &=& 0 ,\qquad  (t, x )
\in (0 , \infty ) \times\mathbb{R}
\nonumber
\\[-8pt]\\[-8pt]
\bar{u} ( 0 ,x ) &=& - \psi ( x ),\qquad  x \in\mathbb{R}.\nonumber
\end{eqnarray}
Since $\bar{u}$ is a classical solution of this equation on $ [h,\infty ) \times\mathbb{R}$,
$\bar{u}  ( \cdot, \bar{x}+ \cdot )$ also classically satisfies
\[
\partial_t \bar{u} (\cdot, \bar{x}+ \cdot ) - b_{k_{\pm}}
\partial_x \bar{u} (\cdot, \bar{x}+ \cdot ) - \int
_{\mathbb{R}} \tilde{\delta}_z \bar {u} (\cdot, \bar{x}+
\cdot ) F_{k_{\pm}} ( dz ) = 0
\]
on $ [h,\infty ) \times\mathbb{R}$. Then
\[
\hat{u}_{ h_0} ( t, x ) := \frac{\bar{u}
 (t , x+ h_0  ) - \bar{u}  (t , x )}{\llvert  h_0\rrvert   }
\]
is a classical solution of (\ref{uPIDElin}) on $ [h,\infty
 ) \times\mathbb{R}$ as well.

Lemma~\ref{nnlem53} implies that
\[
\bigl\llVert \hat{u}_{ h_0} ( \bar{t} + \cdot, \cdot ) \bigr\rrVert
_{L^1  (  (-1, 0  ] \mapsto L^1  ( \nu
 )  )} \leq \operatorname{Lip} ( \psi ) \int_{\mathbb{R}}
\min \bigl( 1 , \llvert z \rrvert ^{-1-\alpha} \bigr)\,dz
\]
for $\bar{t} > 1$. By Lemma~\ref{powerest}, it follows that for some $B$,
\[
\bigl\llvert \partial_x \hat{u}_{h_0} ( t_0,
x_0 ) \bigr\rrvert + \frac{ \llvert   \partial_x \hat{u}_{h_0}  ( t_0, x_0  ) -
\partial_x \hat{u}_{h_0}  (t_1, x_1  )\rrvert  }{ (
\llvert   t_0 - t_1 \rrvert  ^{1/\alpha} + \llvert   x_0 - x_1 \rrvert
 )^{\gamma} } \leq B
\]
for every $ ( t_0, x_0  )$, $ (t_1, x_1  ) \in
\mathfrak{C}_{1/2 , 1/2}  ( \bar{t} , \bar{x}
 )$ with $\bar{t} > 1$.
Rewriting this in terms of $\bar{u}$, we see that we have found a
$\gamma$-H{\"o}lder estimate for
\[
\frac{ \partial_x \bar{u}  (t_0 , x + h_0
) - \partial_x\bar{u}  (t_0 , x )}{\llvert  h_0\rrvert  }.
\]
By Lemma~\ref{lem56cc}, $\partial_{xx}^2 \bar{u}$ exists and is
bounded on $ ( 1/2, \infty ) \times\mathbb{R}$. By
rescaling, we get that this actually holds on $ [ h, h+1  ]
\times\mathbb{R}$.
\end{appendix}

\section*{Acknowledgements}

We gratefully acknowledge Dennis Kriventsov of the University of Texas
at Austin and Luis Silvestre of the University of Chicago for valuable
conversations regarding interior regularity estimates for PIDEs. We
also thank the anonymous referees for their helpful advice on improving
the paper. This work is supported by the National Science Foundation
under Grant DMS-0955463.




%

\printhistory
\end{document}